\documentclass[10pt]{amsart}
\usepackage[pagebackref]{hyperref}

\renewcommand{\Re}{\operatorname{Re}}
\renewcommand{\Im}{\operatorname{Im}}
\renewcommand{\d}{\mathrm{d}}
\renewcommand{\L}{\Lambda}

\newcommand{\ts}{\textstyle }

\newcommand{\sLag}{special Lagrangian}

\newcommand{\cH}{{\mathcal H}}
\newcommand{\cI}{{\mathcal I}}
\newcommand{\cL}{{\mathcal L}}
\newcommand{\cS}{{\mathcal S}}

\newcommand{\bbR}{{\mathbb R}}
\newcommand{\R}[1]{\ensuremath{\mathbb R^{\,#1}}} 
\newcommand{\bbC}{{\mathbb C}}
\newcommand{\C}[1]{\ensuremath{\mathbb C^{\,#1}}} 
\newcommand{\bbZ}{{\mathbb Z}}
\newcommand{\bbP}{{\mathbb P}}

\newcommand{\iC}{\mathrm{i}}
\newcommand{\eul}{\mathrm{e}}
\newcommand{\ov}{\overline}
\newcommand{\ot}{\otimes}

\newcommand{\eusu}{\operatorname{\mathfrak{su}}}
\newcommand{\euu}{\operatorname{\mathfrak u}}

\newcommand{\GL}{\operatorname{GL}}
\newcommand{\SO}{\operatorname{SO}}
\newcommand{\SU}{\operatorname{SU}}

\newcommand{\G}{\operatorname{G}}

\DeclareMathOperator{\tr}{tr}

\newcommand{\bb}{\mathbf{b}}
\newcommand{\eb}{\mathbf{e}}
\newcommand{\ub}{\mathbf{u}}
\newcommand{\vb}{\mathbf{v}}
\newcommand{\xb}{\mathbf{x}}

\newcommand{\w}{{\mathchoice{\,{\scriptstyle\wedge}\,}{{\scriptstyle\wedge}}
      {{\scriptscriptstyle\wedge}}{{\scriptscriptstyle\wedge}}}}
\newcommand{\lhk}{\mathbin{\hbox{\vrule height1.4pt width4pt depth-1pt 
             \vrule height4pt width0.4pt depth-1pt}}}

\numberwithin{equation}{section}

\newtheorem{theorem}{Theorem}

\newtheorem{proposition}{Proposition}
\newtheorem{corollary}{Corollary}

\theoremstyle{remark}

\newtheorem{remark}{Remark}
\newtheorem{example}{Example}

\begin{document}

\author[R. Bryant]{Robert L. Bryant}
\address{Duke University Mathematics Department\\
         P.O. Box 90320\\
         Durham, NC 27708-0320}
\email{\href{mailto:bryant@math.duke.edu}{bryant@math.duke.edu}}
\urladdr{\href{http://www.math.duke.edu/~bryant}%
         {http://www.math.duke.edu/\lower3pt\hbox{\symbol{'176}}bryant}}

\title[Special Lagrangian $3$-folds]
      {Second order families \\
          of special Lagrangian $3$-folds}

\date{May 1, 2001}

\begin{abstract}
A second order family of \sLag{} submanifolds
of~$\C{m}$ is a family characterized by the satisfaction 
of a set of pointwise conditions on the second fundamental form.
For example, the set of ruled \sLag{} submanifolds
of~$\C{3}$ is characterized by a single algebraic equation on
the second fundamental form.

While the `generic' set of such conditions turns out to
be incompatible, i.e., there are no \sLag{} 
submanifolds that satisfy them, there are many
interesting sets of conditions for which the corresponding
family is unexpectedly large.  In some cases, these geometrically
defined families can be described explicitly, leading to
new examples of \sLag{} submanifolds.  In other cases,
these conditions characterize already known families in a new
way.  For example, the examples of Lawlor-Harvey constructed
for the solution of the angle conjecture and recently generalized
by Joyce turn out to be a natural and easily
described second order family.
\end{abstract}



\thanks{
Thanks to Duke University for its support via a research grant 
and to the National Science Foundation 
for its support via DMS-9870164.\hfill\break
\hspace*{\parindent} 
This is Version~$2$. The most recent version
can be found at arXiv:math.DG/0007128.
}

\maketitle

\setcounter{tocdepth}{2}
\tableofcontents

\section{Introduction}\label{sec: intro}

The study of \sLag{} submanifolds was introduced
by Harvey and Lawson in \S III of their fundamental 
paper~\cite{MR85i:53058} on calibrated geometries.  They
analyzed the local and global geometry of these submanifolds
in flat complex $m$-space, constructing many interesting examples
and proving local existence theorems.

Several important classes of examples of these submanifolds have been 
constructed since then, mainly with an eye to applications in
the theory of calibrations or minimizing submanifolds.  Particularly
important was the construction by Lawlor~\cite{MR95f:49057} of 
a \sLag{} manifold asymptotic to a pair of planes
that violate the angle criterion, thus proving that such a pair
of planes is not area-minimizing.  
See also Harvey~\cite{MR91e:53056} for a thorough account of this
example and its applications.

The deformation theory of compact examples in Calabi-Yau
manifolds was studied in the 1990 thesis of R. McLean~\cite{MR99j:53083}, 
who showed that the moduli space of compact \sLag{} 
submanifolds of a given Calabi-Yau manifold is always a disjoint
union of smooth manifolds.

Special Lagrangian geometry received renewed
attention in~1996 when its role in mirror symmetry was discovered
by Strominger, Yau, and Zaslow~\cite{MR97j:32022}.  
Since then, interest in \sLag{} geometry has grown
quite rapidly.   The reader might consult
\cite{MR2000c:32075} (for what is currently known about
the moduli of compact \sLag{} submanifolds), 
\cite{MR2000f:32033} (for some examples arising from algebraic
geometry), \cite{math.DG/0005118}
(for further information about mirror symmetry),
\cite{math.DG/0003220} (for examples with large symmetry groups), 
\cite{math.DG/9912246} (for embedding a given real-analytic Riemannian
$3$-manifold as a \sLag{} submanifold of a Calabi-Yau
$3$-fold), \cite{hep-th/9907013} (for some interesting speculations
about how one might count the isolated \sLag{} submanifolds),
and \cite{math.DG/0005164} (for information about \sLag{}
cones in~$\C{3}$).

Still, the systematic exploration of \sLag{} geometry
seems to have hardly begun.  The known explicit examples have largely 
been found by the well-known Ansatz of symmetry reduction or other
special tricks.  

The research that lead to this article 
was an attempt to classify families of \sLag{} submanifolds
that are characterized by invariant, differential geometric 
conditions, in particular, conditions 
on the second fundamental form of the \sLag{} submanifold.

At least when the ambient space is flat, 
the lowest order invariant of a \sLag{} 
submanifold is its second fundamental form.  Now, for a
Lagrangian submanifold of a linear symplectic vector space, 
the second fundamental form, usually defined as a quadratic form
with values in the normal bundle, has a natural interpretation 
as a symmetric cubic form~$C$ on the submanifold, called the
\emph{fundamental cubic}.  
When the submanifold is \sLag, it turns out that the trace
of this cubic form with respect to the first fundamental form
vanishes, but there are no further pointwise conditions on
this cubic that are satisfied for all \sLag{} submanifolds.

It is natural to ask whether one can obtain nontrivial families 
of \sLag{} submanifolds by imposing pointwise conditions
on the fundamental cubic.  In the language of overdetermined
systems of PDE, one would like to be able to say whether there
are any second order systems of PDE that are `compatible'
with the (first order) system that represents the \sLag{}
condition.%
\footnote{ Here, `compatibility' is not strictly defined, 
but, roughly speaking, means that there exist at least
as many (local) solutions
to the overdetermined system as one would expect 
from a na{\"\i}ve `equation counting' argument.  A more 
precise description would involve concepts from exterior
differential systems, such as involutivity, that will not 
be needed in this article.}

The first task is to understand the space of pointwise
invariants of a traceless cubic form under the special orthogonal
group.  For example, in dimension~$3$ (which is the case this
article mainly considers), the space of traceless cubics
is an irreducible~$\SO(3)$-module of dimension~$7$, so one
would expect there to be four independent polynomial
invariants.  

However, the relations that one gets by 
imposing conditions on these invariants are generally 
singular at the cubics that have a nontrivial stabilizer
under the action of~$\SO(3)$.  For comparison, consider
the classical case of hypersurfaces in Euclidean space.
The fundamental invariants are the principal curvatures,
i.e., the eigenvalues of the second fundamental form 
with respect to the first fundamental form.  These are smooth
away from the (generalized) umbilic locus, i.e., the places
where two or more of the principal curvatures come together.
It is exactly at these places that the stabilizer of the
second fundamental form in the orthogonal group is larger
than the minimum possible stabilizer.  Of course,
the umbilic locus is also the place where moving frame adaptations
generally run into trouble, unless one assumes that the 
multiplicities of the principal curvatures are constant.

There is a similar phenomenon in \sLag{} geometry.  In place
of the umbilic locus, one looks that the places where
the fundamental cubic has a nontrivial stabilizer,%
\footnote{ In contrast to the familiar case of hypersurfaces
in Euclidean space, where the stabilizer, though generically
finite, is always
nontrivial, it turns out that the stabilizer at a 
generic point of the fundamental cubic of a 
`generic' \sLag{} $3$-fold is trivial.} and at the
\sLag{} submanifolds where the stabilizer of the cubic
is nontrivial at the generic point.  These are the special
\sLag{} submanifolds.  

In this article, after making some general remarks to
introduce the structure equations of \sLag{} geometry,
I classify the possible nontrivial $\SO(3)$-stabilizers 
of traceless cubics in three variables.  It turns out that
the $\SO(3)$-stabilizer of a nontrivial traceless cubic 
is isomorphic to either a copy of~$\SO(2)$, 
the group~$A_4$ of order~$12$, the group~$S_3$ of order~$6$, 
the group~$\bbZ_3$, the group~$\bbZ_2$, or is trivial.
I then consider, for each of the nontrivial subgroups~$G$
on this list, the problem of classifying the \sLag{} $3$-folds 
whose cubic form at each point has its stabilizer contain a
copy of~$G$.

For example, it turns out that the only \sLag{} $3$-folds
in~$\C3$ whose cubic form has a continuous stabilizer at
each point are the $3$-planes and the $\SO(3)$-invariant
examples discovered by Harvey and Lawson.

There are no \sLag{} $3$-folds whose cubic stabilizer at
at generic point is of type~$A_4$, but the ones whose stabilizer
at a generic point is of type~$S_3$ turn out to be the austere
\sLag{} $3$-folds and these are known to be the orthogonal
products, the \sLag{} cones, and the `twisted' \sLag{} cones.

The \sLag{} $3$-folds with cubic stabilizer at a generic
point isomorphic to~$\bbZ_2$ turn out to be the examples
discovered by Lawlor, extended by the work of Harvey and
then Joyce.%
\footnote{ I am indebted to Joyce for suggesting 
(by private communication) that the
family that I had shown to exist in this case might be
the Lawlor-Harvey-Joyce family.  He was correct, and this
saved me quite a bit of work in integrating the corresponding
structure equations.}

The \sLag{} $3$-folds whose cubic stabilizer at a generic
point is isomorphic to~$\bbZ_3$ turn out to be asymptotically
conical and, indeed, turn out to be deformations, in a certain
sense, of the \sLag{} cones, as explained in the 
thesis of Haskins~\cite{math.DG/0005164}.

The above results are explained more fully 
in~\S\ref{sec: 2nd order fam}.  

At the conclusion of \S\ref{sec: 2nd order fam}, 
I consider a different type of invariant condition
on the fundamental cubic, namely that, at every point~$x\in L$,
the degree~$3$ curve in~$\bbP(T_xL)$ defined by the fundamental
cubic have a real singular point.  This is one semi-algebraic 
condition on the fundamental cubic.
I show that the \sLag{} $3$-folds with this property are
exactly the ruled \sLag{} $3$-folds.  

Moreover, while 
the most general known family of ruled \sLag{} $3$-folds
up until now was one discovered by Borisenko~\cite{MR94f:53099} and
that depends on four functions of one variable (in the
sense of exterior differential systems), I show that the
full family depends on six functions of one variable.
Moreover, I show that, when one interprets the ruled
\sLag{} $3$-folds as surfaces in the space of lines in~$\C3$,
the surfaces that one obtains are simply the ones that
are holomorphic with respect to a canonical Levi-flat,
almost CR-structure on the space of lines.  This interpretation
has several implications for the structure of ruled \sLag{}
$3$-folds, among them being that any real-analytic ruled
surface in~$\C3$ on which the K\"ahler form vanishes lies
in a (essentially unique) ruled \sLag{} $3$-fold.  Moreover,
a \sLag{} $3$-fold is ruled if and only if it contains
a ruled surface.

It has to be said that the results of this article are 
only the first step in understanding the compatibility
of the \sLag{} condition with higher order conditions.  
Now that the `umbilic' cases are understood, the serious
work on the `generic' case can be undertaken.  This will
be reported on in a subsequent work.

Also, while, for the sake of brevity, this work has concerned itself 
(essentially exclusively) with the $3$-dimensional case, 
there are obvious higher dimensional generalizations 
that need to be investigated and that should yield to the 
same or similar techniques.%
\footnote{ My student, Marianty Ionel, has recently
completed a study of the \sLag{} $4$-folds in~$\C4$ whose
fundamental cubic has nontrivial symmetries.}

\subsection{Special Lagrangian geometry}\label{ssec: sLag geom}

In this article, a slightly more general notion of \sLag{}
geometry is adopted than is customary.  The reader might compare
this discussion with Harvey and Lawson's original article~\cite{MR85i:53058} 
or Harvey's more recent book~\cite{MR91e:53056}.

\subsubsection{Special K\"ahler structures}\label{sssec: sK strs}

Let $M$ be a complex $m$-manifold endowed with a K\"ahler form~$\omega$
and a holomorphic volume form~$\Upsilon$.  It is \emph{not} assumed
that~$\Upsilon$ be parallel, or even of constant norm, with respect
to the Levi-Civita connection associated to~$\omega$.  
The pair~$(\omega,\Upsilon)$ is said to define a
\emph{special K\"ahler structure} on~$M$.

\subsubsection{Special Lagrangian submanifolds}\label{sssec: sLags}

A submanifold~$L\subset M$ of real dimension~$m$ is said to be
\emph{Lagrangian}%
\footnote{ or \emph{$\omega$-Lagrangian} 
if there is any danger of confusion} 
if the pullback of~$\omega$ to~$L$ vanishes.

Harvey and Lawson show~\cite[\S III, Theorem~1.7]{MR85i:53058} 
that for any Lagrangian
submanifold~$L\subset M$, the pullback of~$\Upsilon$ to~$L$
can never vanish.  A Lagrangian
submanifold~$L$ is said to be \emph{special Lagrangian} 
if the pullback of~$\Im(\Upsilon)$ to~$L$ vanishes.  
When~$L$ is \sLag{}, it has
a canonical orientation for which~$\Re(\Upsilon)$ pulls back to~$L$
to be a positive volume form, and this is the orientation that 
will be assumed throughout this article.

More generally, if~$\lambda$ is a complex number of unit
modulus, one says that an oriented Lagrangian submanifold~$L\subset M$
has \emph{constant phase}~$\lambda$ if~${\bar\lambda}\,\Upsilon$ 
pulls back to~$L$ to be a (real-valued) positive volume form.
Obviously, for any fixed~$\lambda$, this notion is not significantly
more general than the notion of \sLag{}, so I will
usually consider only \sLag{} submanifolds in this article.

\subsubsection{The Calabi-Yau case}\label{sssec: Calabi-Yau}

When~$\Upsilon$ is parallel with respect to the Levi-Civita
connection associated to~$\omega$, the K\"ahler metric has
vanishing Ricci tensor and the pair~$(\omega,\Upsilon)$ is
said to define a \emph{Calabi-Yau} structure on~$M$.  In this
case, Harvey and Lawson show that any \sLag{} 
submanifold~$L\subset M$ is minimal.  Moreover, if~$L$ is
compact, it is absolutely minimizing in its homology 
class since it is then calibrated by~$\Re(\Upsilon)$.

\subsubsection{Local existence}\label{sssec: local existence}

Assume that~$\omega$ is real-analytic with respect to the
standard real-analytic structure on~$M$ that underlies its
complex analytic structure.
Harvey and Lawson~show\cite[\S III, Theorem~5.5]{MR85i:53058} 
that any real-analytic submanifold~$N\subset M$
of dimension~$m{-}1$ on which~$\omega$ pulls back to be zero 
lies in a unique \sLag{} submanifold~$L\subset M$.  
(Although their result is stated only for the case of the standard
flat special K\"ahler structure on~$\C{m}$, their proof 
is valid in the general case, provided one makes the necessary
trivial notational changes.)

Thus, there are many \sLag{} submanifolds locally,
at least in the real-analytic category.  By adapting arguments
from~\cite[\S III.2]{MR85i:53058}, one can also prove local
existence of \sLag{} submanifolds even without
the assumption of real-analyticity.  Instead, one uses local
existence for an elliptic second order scalar equation.

\subsubsection{Deformations}\label{sssec: deformations}

R. McLean~\cite{MR99j:53083} proved that, 
in the Calabi-Yau case, a compact \sLag{}
submanifold~$L\subset M$ is a point in a smooth, finite
dimensional moduli space~$\cL$ consisting of the 
\sLag{} deformations of~$L$ and that the tangent space
to~$\cL$ at~$L$ is isomorphic to the space of harmonic
1-forms on~$L$.  

McLean's argument makes no essential use of the assumption 
that~$\Upsilon$ be $\omega$-parallel.
Instead, it is sufficient for the conclusion of 
McLean's theorem that~$\Upsilon$ be closed (in fact, one only
really needs that the imaginary part of~$\Upsilon$ be closed.)
For a related result, see~\cite{math.DG/9906048}.

\subsection{Special K\"ahler reduction}\label{ssec: sK reduction}

One reason for considering the slightly wider notion of 
\sLag{} geometry adopted here is that it is
stable under the process of \emph{reduction}, as explained
in~\cite{math.DG/0002097}, \cite{math.DG/0008021}, 
and~\cite{math.AG/0012002}.

Let~$(\omega,\Upsilon)$ be a special K\"ahler structure on~$M$.
A vector field~$X$ on~$M$ will be said to be an 
\emph{infinitesimal symmetry} of the structure 
if the (locally defined) flow of~$X$ preserves both~$\omega$
and~$\Upsilon$.  

Suppose that~$X$ is an infinitesimal symmetry of~$(\omega,\Upsilon)$
and that $X$ is, moreover, $\omega$-Hamiltonian, 
i.e., that there exists a function~$H$ on~$M$ 
satisfying~$X\lhk\omega = -\d H$.  The flow lines of~$X$ are
tangent to the level sets of~$H$.  

Say that a value~$h\in\bbR$
is a \emph{good} value for~$H$ if it is a regular value of~$H$
and if the flow of~$X$ on the 
level set~$H^{-1}(h)\subset M$ is simple, 
i.e., there is a smooth manifold structure on the set~$M_{h}$ 
of flow lines of~$X$ in the level set~$H^{-1}(h)$ so that the
natural projection~$\pi_h:H^{-1}(h)\to M_h$ is a smooth submersion.
The (real) dimension of~$M_h$ is necessarily~$2m{-}2$.

When~$h$ is good, there exists a unique $2$-form~$\omega_h$ 
on~$M_h$ for which~$\pi_h^*(\omega_h)$ is the pullback of~$\omega$
to~$H^{-1}(h)$ and there exists a unique complex-valued~$(m{-}1)$
form~$\Upsilon_h$ on~$M_h$ for which~$\pi_h^*(\Upsilon_h)$ is
the pullback to~$H^{-1}(h)$ of~$X\lhk\Upsilon$.  

It is trivial to verify that~$(\omega_h,\Upsilon_h)$ defines
a special K\"ahler structure on~$M_h$.  Note, however, that, even
if~$(\omega,\Upsilon)$ is Calabi-Yau, its reductions will generally
\emph{not} be Calabi-Yau.  In fact, this happens only when the
length of~$X$ is constant along the level set~$H^{-1}(h)$.

If~$L\subset H^{-1}(h)$ is a \sLag{} submanifold that
is tangent to the flow of~$X$, then~$L = \pi_h^{-1}(L_h)$
where~$L_h\subset M_h$ is also \sLag{}.  Conversely,
if $L_h\subset M_h$ is \sLag{}, then $L = \pi_h^{-1}(L_h)$
is \sLag{} in~$M$. 

This method of special K\"ahler reduction allows one to construct
many examples of \sLag{} submanifolds by starting
with a Hamiltonian $(m{-}1)$-torus action and doing a 
series of reductions, leading to a 1-dimensional 
special K\"ahler manifold, where the integration problem is 
reduced to integrating a holomorphic 1-form on a Riemann surface.

\section[Structure Equations]{The Structure Equations}\label{sec: str eqs}

The structure equations of a K\"ahler manifold adapted for 
\sLag{} geometry can be found in~\cite{MR88j:53061},
and will only be reviewed briefly here.

\subsection{The special coframe bundle}\label{ssec: spec cof bundl}

The \emph{standard special K\"ahler structure} on~$\C{m}$ is
the one defined by 
\begin{equation}
\omega_0 = {\ts\frac\iC2}\bigl(\d z_1\w\d \ov{z_1} 
             + \dots + \d z_m\w\d\ov{z_m}\bigr)
\qquad\text{and}\qquad
\Upsilon_0 = \d z_1\w\dots\w\d z_m\,,
\end{equation}
where~$z_1,\dots,z_m$ are the usual complex linear coordinates on~$\C{m}$.
The corresponding K\"ahler metric is, of course
\begin{equation}
g_0 = \d z_1\circ\d \ov{z_1} 
             + \dots + \d z_m\circ\d \ov{z_m}\,.
\end{equation}
Note that~$\R{m}\subset\C{m}$ is a \sLag{} subspace.

Let~$(\omega,\Upsilon)$ be a special K\"ahler structure on
the complex $m$-manifold~$M$.  There is a unique positive function~$B$
on~$M$ that satisfies
\begin{equation}
\Upsilon\w\ov{\Upsilon} = \frac{2^m (-\iC)^{m^2}}{m!}\,B^2\,\omega^m\,.
\end{equation}

A linear isomorphism~$u:T_xM\to\C{m}$ will be said to be a
\emph{special K\"ahler coframe} at~$x$ if it 
satisfies~$\omega_x=u^*(\omega_0)$ 
and~$\Upsilon_x = B(x)\,u^*(\Upsilon_0)$.  Such a coframe
is necessarily complex linear, i.e., satisfies~$u(J_xv) = \iC\,u(v)$
for all~$v\in T_xM$, where~$J_x:T_xM\to T_xM$ is the complex structure
map.

The set of special K\"ahler coframes at~$x$ will be denoted~$P_x$
and is the fiber of a principal right~$\SU(m)$-bundle~$\pi:P\to M$,
with right action given by~$R_a(u) = a^{-1}\circ u$ for
$a\in\SU(m)$.

As usual, the $\C{m}$-valued, tautological 1-form~$\zeta$ on~$P$ 
is defined by requiring that~$\zeta_u = u\circ\pi'(u):T_uP\to \C{m}$
for~$u\in P$.  It satisfies~$R_a^*(\zeta) = a^{-1}\,\zeta$ 
for~$a\in\SU(m)$.

The components of~$\zeta$ will be written as~$\zeta_i$ for~$1\le i\le m$. 
The equations
\begin{equation}
\omega = {\ts\frac\iC2}\bigl(\zeta_1\w\ov{\zeta_1} 
             + \dots + \zeta_m\w\ov{\zeta_m}\bigr)
\qquad\text{and}\qquad
\Upsilon = B\,\zeta_1\w\dots\w\zeta_m
\end{equation}
hold on~$P$, where, as is customary, I have omitted the $\pi^*$,
thus implicitly embedding the differential forms on~$M$ into the
differential forms on~$P$ via pullback.

Finally, there are functions~$\eb = (\eb_i)$, where~$\eb_i:P\to TM$
is a bundle mapping satisfying~$\zeta_i(\eb_j) = \delta_{ij}$.
In other words~$\pi'(v) = \zeta_i(v)\,\eb_i(u)$ for all~$v\in T_uP$.

\begin{remark}[The flat case]
When~$M=\C{m}$ and $(\omega,\Upsilon) = (\omega_0,\Upsilon_0)$, it
is customary to use the vector space (parallel) trivialization 
of the tangent bundle of~$\C{m}$ to identify all of the tangent spaces
to the vector space~$\C{m}$ itself.  In this case, the functions~$\eb_i$
will be regarded as vector-valued functions on~$P\simeq\C{m}\times\SU(m)$
and the basepoint projection will be denoted as~$\xb:P\to\C{m}$.  
Then the above relations take on the more familiar `moving frame'
form
$$
\d\xb = \eb_i\,\zeta_i\,,
$$
and so on.  The reader should have no trouble figuring out what is
meant in context.
\end{remark}

\subsection{The structure equations}\label{ssec: cmplex str eqs} 

The Levi-Civita connection associated to the underlying K\"ahler
structure on~$M$ is represented on~$P$ by a $\euu(m)$-valued 1-form
$\psi = -\psi^* = (\psi_{i\bar\jmath})$ that satisfies the 
\emph{first structure equation}
\begin{equation}
\d\zeta_i = -\psi_{i\bar\jmath}\w\zeta_{j}\,.
\end{equation}
The equation~$\d\Upsilon=0$ implies
\begin{equation}\label{eq: vol torsion}
(\bar\partial - \partial)\bigl(\log B\bigr)
= -\iC\,\d^c(\log B) = \tr(\psi) = \psi_{i\imath}\,.
\end{equation}
Note that~$(\omega,\Upsilon)$ is Calabi-Yau if and only if~$B$
is constant, i.e., if and only if~$\psi$ takes values in~$\eusu(m)$.

In the Calabi-Yau case, where~$\Upsilon$ is parallel with respect to
the Levi-Civita connection of~$\omega$, the relation~$\psi_{ii}=0$
holds.  Moreover, the Calabi-Yau structure is locally equivalent
to the standard structure if and only if the Levi-Civita connection
of~$\omega$ vanishes, i.e., if and only if $\d\psi = -\psi\w\psi$,
which is known as the \emph{second structure equation} of a flat 
Calabi-Yau space.

\subsection{Special Lagrangian submanifolds}\label{ssec: str of slag}

For the study of \sLag{} submanifolds, it is convenient
to separate the structure equations into real and imaginary parts.
Thus, set~$\zeta_i = \omega_i + \iC\,\eta_i$ and~$\psi_{i\bar\jmath}
= \alpha_{ij} + \iC\,\beta_{ij}$.  The first structure equations
can then be written in the form
\begin{equation}\label{eq: real form str eqs}
\begin{split}
\d\,\omega_i &= - \alpha_{ij}\w\omega_j +  \beta_{ij}\w\eta_j\,,\\
\d\,\eta_i   &= -  \beta_{ij}\w\omega_j - \alpha_{ij}\w\eta_j\,.
\end{split}
\end{equation}
where~$\alpha_{ij} = -\alpha_{ji}$ and~$\beta_{ij} = \beta_{ji}$.  
Note that~\eqref{eq: vol torsion} becomes~$\beta_{ii} 
= -\d^c\bigl(\log B\bigr)$.

Let~$L\subset M$ be a \sLag{} submanifold.  For~$x\in L$, 
a special K\"ahler coframe~$u:T_xM\to\C{m}$ is said to 
be \emph{$L$-adapted} if~$u(T_xL)=\R{m}\subset\C{m}$
and~$u:T_xL\to\R{m}$ is orientation preserving. 
The space of $L$-adapted coframes forms
a principal right~$\SO(m)$-subbundle~$P_L\subset\pi^{-1}(L)\subset P$
over~$L$.

The equations~$\eta_i=0$ hold on~$P_L$.  Thus, by the
structure equations~\eqref{eq: real form str eqs}, 
the relations~$\beta_{ij}\w\omega_j = 0$ hold on~$P_L$ while~$\omega_1\w\dots
\w\omega_m$ is nowhere vanishing.  
It follows from Cartan's Lemma that there are functions~$h_{ijk} = h_{jik}
= h_{ikj}$ on~$P_L$ so that
\begin{equation}
\beta_{ij} = h_{ijk}\,\omega_k\,.
\end{equation}
The second fundamental form of~$L$ can then be written as
\begin{equation}
\mathrm{I\!I}= h_{ijk}\,J\eb_i\,{\ot}\,\omega_j\omega_k 
               = J\eb_i\,{\ot}\,Q_i\,
\end{equation}
where~$Q_i = h_{ijk}\,\omega_j\omega_k$.
The information in the second fundamental form
is thus contained in the symmetric cubic form
\begin{equation}
C = h_{ijk}\,\omega_i\omega_j\omega_k = \omega_i\,Q_i\,,
\end{equation}
which is well-defined on~$L$.  This symmetric cubic form will be
referred to as the \emph{fundamental cubic} 
of the \sLag{} submanifold~$L$.

Note that the trace of~$C$
with respect to the 
induced metric~$g = {\omega_1}^2 +\dots+{\omega_n}^2$ on~$L$ satisfies
\begin{equation}
\tr_gC = h_{iik}\,\omega_k = \beta_{ii} 
= -\d^c(\log B)\vrule\,{}_L\,,
\end{equation}
which is the restriction to~$L$ of an ambient~$1$-form.
In the Calabi-Yau case, $0 = \psi_{ii} = \iC\,\beta_{ii}$, so the
fundamental cubic~$C$ is traceless.  

Finally, in the flat case, the curvature vanishing 
condition~$\d\psi = -\psi\w\psi$ can be separated into real and imaginary
parts.  The result will be referred to 
as the \emph{second structure equations}:
\begin{subequations}
\begin{align}
\d\alpha_{ij} &= - \alpha_{ik}\w\alpha_{kj}+  \beta_{ik}\w\beta_{kj}\,,
\label{eq: 2nd str eq 1 -- flat}\\
\d\beta_{ij}  &= -  \beta_{ik}\w\alpha_{kj}- \alpha_{ik}\w\beta_{kj}\,.
\label{eq: 2nd str eq 2 -- flat}
\end{align}
\end{subequations}

\subsubsection{A Bonnet-type result}\label{ssec: slag bonnet}
Given an $m$-manifold~$L$ endowed with a Riemannian metric~$g$ and
a symmetric cubic form~$C$ that is traceless with respect to~$g$, 
one can ask whether there is an isometric imbedding of~$(L,g)$ 
into~$\C{m}$ as a \sLag{} submanifold that induces~$C$ as the 
fundamental cubic.  

It is easy to see that there is a Bonnet-style
theorem derivable from the above structure equations.  Namely,
there are analogs of the Gauss and Codazzi equations that give
necessary and sufficient conditions for the solution of this
problem.

To see this, first choose a $g$-orthonormal coframing~$\omega = (\omega_i)$
on an open subset~$U\subset L$.  Then define~$\eta_i=0$ and let~$\alpha_{ij}
=-\alpha_{ji}$ be the unique $1$-forms on~$U$ that satisfy the
equations~$\d\omega_i = -\alpha_{ij}\w\omega_j$.  (The existence and 
uniqueness of such~$\alpha$ is just the Fundamental Lemma of Riemannian
geometry.) 
Write~$C = h_{ijk}\,\omega_i\omega_j\omega_k$ and 
set~$\beta_{ij} = h_{ijk}\,\omega_k$.  

The equation~\eqref{eq: 2nd str eq 2 -- flat} then expresses
the fact that~$C$ must satisfy a Codazzi-type equation with respect to~$g$,
namely, that the covariant derivative of~$C$ with respect to the Levi-Civita
connection of~$g$ is fully symmetric. 
The equation~\eqref{eq: 2nd str eq 1 -- flat} then
expresses the fact that that~$C$ must satisfy a Gauss-type 
equation with respect to~$g$, namely, that~$C$ satisfies an
algebraic equation of the form~$Q_g(C)=\text{Riem}(g)$, where~$Q_g$
is a certain quadratic mapping (depending on~$g$) 
from symmetric cubic forms into tensors of the same 
algebraic type as the Riemann curvature.

Thus, when~$L$ is simply connected, these Codazzi and Gauss equations 
are the necessary and sufficient conditions on~$g$ and~$C$ 
for there to be an isometric immersion of~$(L,g)$ into~$\C{m}$
as a \sLag{} submanifold inducing~$C$ as its fundamental cubic.  
Moreover, such an isometric immersion will be unique up to rigid motion.

\section{Second Order Families}\label{sec: 2nd order fam}

\subsection{The second fundamental form as a cubic}\label{ssec: cubic}

It was already explained in~\S\ref{sec: str eqs} how the second
fundamental form of a \sLag{} submanifold~$L\subset\C{m}$ can be
regarded as a symmetric cubic form that is traceless with respect
to the first fundamental form.  

Thus, the polynomial second order 
invariants of such a submanifold correspond to the $\SO(m)$-invariant
polynomials on the space~$\cH_3(\R{m})$ 
of harmonic polynomials on~$\R{m}$ that are homogeneous of degree~$3$.
Moreover, the $\SO(m)$-stabilizer of a given cubic in this space
corresponds to the ambiguity in the choice of an adapted coframe
for the corresponding \sLag{} submanifold.  In particular, points
on \sLag{} submanifolds
at which the $\SO(m)$-stabilizer of the second fundamental form is 
nontrivial can be regarded as analogs of umbilic points in the
classical theory of surfaces in 3-space.

The space~$\cH_3(\R{m})$ is an irreducible~$\SO(m)$-module when $\SO(m)$
acts in the natural way by pullback.  Thus, there are no invariant
linear functions and, up to multiples, exactly one invariant quadratic
polynomial, which is, essentially, the squared norm of the polynomial.
It is not difficult to show that there are no invariant cubic polynomials
on this space, that the space of invariant quartics is of dimension~$2$
(one of which is the square of the invariant quadratic), and so on.
The exact structure of the ring of invariants for general~$m$ 
is complicated, however, and I will not discuss this further.

Since I will only be using the results of the case~$m=3$, I am going
to be assuming this from now on.  The space~$\cH_3(\R{3})$ has
dimension~$7$, and one would expect that the `generic' 
$\SO(3)$-orbit in this vector space to have dimension~$3$.  
In fact, I am now going to determine the orbits that have 
non-trivial stabilizers.

\subsubsection{Special orbits}\label{sssec: special orbits}

The main goal of this section is to prove the following result,
which is undoubtedly classical even though I have been unable to locate
a proof in the literature.

\begin{proposition}  \label{prop: nontrivial stabilizers}
The $\SO(3)$-stabilizer of~$h\in\cH_3(\R{3})$ 
is nontrivial if and only if~$h$ lies on the~$\SO(3)$-orbit of
exactly one of the following polynomials
\begin{enumerate}
\item $0\in\cH_3(\R{3})$, whose stabilizer is~$\SO(3)$.
\item $r\,(2z^3-3zx^2-3zy^2)$ for some~$r>0$, whose stabilizer is~$\SO(2)$.
\item $6s\,xyz$ for some~$s>0$, whose stabilizer is the
      subgroup~$A_4\subset\SO(3)$ of order~$12$ generated by
      the rotations by an angle of~$\pi$ about the $x$-, $y$-, and $z$-axes
      and by rotation by an angle of~$\frac23\pi$ about the line~$x=y=z$.
\item $s\,(x^3-3xy^2)$ for some~$s>0$, whose stabilizer is the
      subgroup~$S_3\subset\SO(3)$ of order~$6$ generated
      by the rotation by an angle of~$\pi$ about the~$x$-axis 
      and the rotation by an angle of~$\frac23\pi$ about the~$z$-axis. 
\item $r\,(2z^3-3zx^2-3zy^2)+6s\,xyz$ for some~$r,s>0$ satisfying~$s\not=r$,
      whose stabilizer is the~$\bbZ_2$-subgroup of~$\SO(3)$ generated by
      rotation by an angle of~$\pi$ about the~$z$-axis.
\item $r\,(2z^3-3zx^2-3zy^2)+s\,(x^3-3xy^2)$ 
      for some~$r,s>0$ satisfying~$s\not=r{\sqrt2}$,
      whose stabilizer is the~$\bbZ_3$-subgroup of~$\SO(3)$ generated by
      rotation by an angle of~$\frac23\pi$ about the~$z$-axis.
\end{enumerate}
\end{proposition}

\begin{remark}[Special Values]\label{rem: special values}
The reader may wonder about the conditions~$s\not=r$ and~$s\not=r\sqrt2$
in the last two cases.  It is not difficult to verify that the 
polynomial $(2z^3-3zx^2-3zy^2)+6xyz$ lies on the
$\SO(3)$-orbit of~$2(x^3-3xy^2)$ and that the polynomial
$(2z^3-3zx^2-3zy^2)+{\sqrt2}\,(x^3-3xy^2)$ lies on the $\SO(3)$-orbit
of~$6{\sqrt3}\,xyz$.
\end{remark}

\begin{proof}
Suppose that~$h\in\cH_3(\R{3})$ has a nontrivial stabilizer~$G\subset\SO(3)$.
Obviously~$G=\SO(3)$ if and only if~$h=0$, so suppose that~$h\not=0$
from now on.    Since~$G$ is closed in~$\SO(3)$, it is compact and has
a finite number of components.  

Suppose first that~$G$ is not discrete.  
Then the identity component of~$G$ must be a closed $1$-dimensional 
subgroup and hence conjugate to the subgroup~$\SO(2)\subset\SO(3)$
consisting of the rotations about the~$z$-axis.  Thus,~$h$ lies on
the orbit of a cubic polynomial that is invariant under this rotation
group.  By replacing~$h$ by such an element, it can be supposed that
the identity component of~$G$ is~$\SO(2)$.  

Consider the following four subspaces of~$\cH_3(\R{3})$:  Let~$V_0$ 
be the $1$-dimensional space spanned by~$z(2z^2{-}3x^2{-}3y^2)$;
let~$V_1$ be the $2$-dimensional space spanned by $x(4z^2{-}x^2{-}y^2)$
and~$y(4z^2{-}x^2{-}y^2)$; let~$V_2$ be the $2$-dimensional space
spanned by~$(x^2{-}y^2)z$ and $xyz$; and let~$V_3$ be the $2$-dimensional
space spanned by $(x^3{-}3xy^2)$ and $(3x^2y{-}y^3)$.  Each of
these subspaces is preserved by the elements of~$\SO(2)$.  Moreover,
$\SO(2)$ acts trivially on~$V_0$, while the element~$R_\alpha\in\SO(2)$
that represents rotation by an angle~$\alpha$ about the $z$-axis, acts
as rotation by the angle~$k\alpha$ on the $2$-dimensional space~$V_k$
for~$k=1,2,3$.  

Obviously, the only nonzero elements of~$\cH_3(\R{3})$ that
are fixed by~$\SO(2)$ are those of the form~$rz(2z^2{-}3x^2{-}3y^2)$
for some nonzero~$r$.  Moreover, since~$z$ is the unique linear
factor of this polynomial, it follows that~$G$, the stabilizer of this
polynomial, must preserve the~$z$-axis.  Also, since this polynomial
is positive on exactly one of the two rays in the $z$-axis emanating
from the origin, it follows that~$G$ must also fix the orientation of 
the~$z$-axis.  Thus, $G=\SO(2)$. Moreover, note that by a rotation
that reverses the~$z$-axis, the element~$rz(2z^2{-}3x^2{-}3y^2)$
is carried into the element~$-rz(2z^2{-}3x^2{-}3y^2)$.  Thus, one
can assume that~$r>0$.

Now suppose that~$G$ is discrete (and hence finite).  Let~$A\in G$
be an element of finite order~$p>1$.  Then~$A$ is rotation about a line
by an angle of the form~$(2q/p)\pi$ for some integer~$q$ relatively
prime to~$p$ and satisfying~$0<q<p$.   Replacing~$h$
by an element in its $\SO(3)$-orbit, I can assume that the fixed 
line of~$A$ is the $z$-axis.  Since the action of~$A$
on~$V_k$ is a rotation by the angle~$(2kq/p)\pi$ for~$k=1,2,3$, it follows
that, unless either~$2q/p$ or~$3q/p$ are integers, then the only 
elements of~$\cH_3(\R{3})$ that are fixed by~$A$ are the elements of~$V_0$.
Since these elements have a continuous symmetry group, and so, by hypothesis,
cannot be~$h$, it follows that either~$2q/p$ or~$3q/p$ are integers, i.e.,
that $p=2$ or $p=3$.

If~$p=2$, then~$h$ must lie in $V_0+V_2$, i.e., 
there must be constants~$r$, $s$, and $t$, so that
$$
h = r\,z(2z^2 - 3x^2 - 3y^2) + 3\bigl(s\,(2xy) + t\,(x^2{-}y^2)\bigr)z\,.
$$
By a rotation that reverses the~$z$-axis, if necessary, I can assume 
that~$r\ge0$ and then, by applying a rotation in~$\SO(2)$, 
I can assume that~$t=0$ and~$s\ge0$.  Since~$G$ is discrete, $s$ cannot
be zero, so~$s>0$.  Note that~$A$ is a rotation by an angle of~$\pi$
about the~$z$ axis, and that this certainly preserves any~$h$ in
the above form.   Note also that every such~$h$ has a linear factor.
In particular, to each element~$A$ of order~$2$ in~$G$, 
there corresponds a linear factor of~$h$ that is fixed (up to a
sign) by~$A$.  

If~$r=0$, then~$h = 6s\,xyz$, and it is clear that
the elements of~$G$ must permute the planes~$x=0$, $y=0$, and~$z=0$.
It follows that~$G$ must be the group~$A_4$ of order~$12$
described in the proposition.  

If~$r>0$, then it is still true that~$h$ has a linear factor, i.e.,
$$
h = \bigl(2r\,z^2 - 3r\,x^2 - 3r\,y^2 + 6s\,xy \bigr)\,z.
$$

When~$r \not= s$, the quadratic factor in the above expression
is irreducible (since $r$ and~$s$ are positive), 
so~$G$ must stabilize the $z$-axis.
In fact, since $h$ is positive on the positive ray of the~$z$-axis,
$G$ must actually be a subgroup of~$\SO(2)$.
Since~$s>0$, $G$ must therefore be isomorphic to~$\bbZ_2$,
 generated by the rotation by~$\pi$ about the~$z$-axis.  
On the other hand, when~$r=s$, the polynomial~$h$ factors as
$$
h =  r({\sqrt2\,z} - {\sqrt3\,x} + {\sqrt3\,y})
      ({\sqrt2\,z} + {\sqrt3\,x} - {\sqrt3\,y})\,z\,.
$$
These three linear factors of~$h$ are linearly dependent, so that
$h$ vanishes on the union of three coaxial planes that meet pairwise
at an angle of~$\pi/3$.  Consequently,~$h$
lies on the $\SO(3)$-orbit of an element of the form
$$
p\,(x^3-3xy^2) = p\,x\,\bigl(x-{\sqrt3\,y}\bigr)\bigl(x+{\sqrt3\,y}\bigr)
$$
where~$p>0$.   Since~$G$ must preserve these factors up to a sign,
$G$ is isomorphic to~$S_3$ and is generated
as claimed in the proposition.

Finally, assume that~$G$ has no element of order~$2$.  Then, 
by the above argument, all of the nontrivial elements of~$G$ have order~$3$.  
By the well-known classification of the finite subgroups of~$\SO(3)$, %
\footnote{ Up to conjugation, these subgroups consist 
of the cyclic subgroups, the dihedral subgroups,
and the symmetry groups of the Platonic solids.} 
it follows that~$G$ must be isomorphic to~$\bbZ_3$.  

Let~$A$ be a generator of~$G$ and assume (as one may, by replacing~$h$
by an element in its~$\SO(3)$-orbit) that~$A$ is rotation by an
angle of~$2\pi/3$ about the~$z$-axis.  Then the elements of~$\cH_3(\R3)$
that are fixed by~$A$ are the elements in~$V_0+V_3$, i.e.,
those  of the form
$$
h = r\,z(2z^2 - 3x^2 - 3y^2) + s\,(x^3-3xy^2) + t\,(3x^2y-y^3).
$$
By a rotation about the~$z$-axis, $h$ can be replaced by an element
in its orbit that is of the above form but that satisfies~$t = 0$
and~$s\ge0$.  Now, $s>0$, since, otherwise the stabilizer of~$h$
would contain~$\SO(2)$.  After rotation by an angle of~$\pi$
about the~$x$-axis if necessary, I can further assume that~$r\ge0$.  
In fact,~$r>0$, since, otherwise,~$G$ would be isomorphic to~$S_3$, 
contrary to hypothesis. 

It remains to determine those 
positive values of~$r$ and~$s$ (if any) for which
$$
h = r\,z(2z^2 - 3x^2 - 3y^2) + s\,(x^3-3xy^2)
$$
has a symmetry group larger than~$\bbZ_3$.  

If the symmetry group~$G$ is to be larger than~$\bbZ_3$, 
then, by the aforementioned classification,
either~$G$ contains an element of order~$2$ or~$G$ is infinite.  In
either case, by the above arguments, $h$ must have a linear factor.
Now, it is straightforward to verify that~$h$ has no linear factor
unless~$s=r\sqrt2$.  Thus, the stabilizer is~$\bbZ_3$ except in
this case.  On the other hand
$$
z(2z^2 - 3x^2 - 3y^2) + {\sqrt2}\,(x^3-3xy^2)
= (z + {\sqrt2}\,x)\bigl(x + {\sqrt3}\,y - {\sqrt2}\,z\bigr)
        \bigl(x - {\sqrt3}\,y - {\sqrt2}\,z\bigr)\,,
$$
and the three linear factors vanish on three mutually orthogonal
$3$-planes.  It follows immediately that this~$h$ lies on the
orbit of~$6p\,xyz$ for some~$p>0$.
\end{proof}

The argument for Proposition~\ref{prop: nontrivial stabilizers}
can be used to prove two more easy results:

\begin{proposition}\label{prop: redcible cubics}
A cubic~$h\in\cH_3(\R3)$ is reducible if and only if it
has a symmetry of order~$2$.  It factors into three linear factors
if and only if it is either the zero cubic, has symmetry~$A_4$,
or has symmetry~$S_3$.
\end{proposition}

\begin{proof}
By Proposition~\ref{prop: nontrivial stabilizers}, any cubic
that has a symmetry of order~$2$ has a linear factor.  Conversely,
suppose that~$h\in\cH_3(\R3)$ has a linear factor and is nonzero. 
By applying an~$\SO(3)$ symmetry, 
it can be assumed that~$z$ divides~$h$, implying that~$h$ has
the form 
$$
h = z\bigl(r(2z^2-3x^2-3y^2) + 3p(x^2-y^2) + 3q(2xy)\bigr),
$$
which clearly has a symmetry of order~$2$ that fixes~$z$.  
The quadratic factor is reducible if and only if either~$r=0$,
in which case a rotation in the $xy$-plane reduces~$p$ to zero,
so that the symmetry group is~$A_4$, 
or else $p^2 + q^2 = r^2$, in which case~$h$ factors into three
linearly dependent factors, so that the symmetry group is~$S_3$.
\end{proof}

Before stating the next proposition, it will be useful to 
establish some notation.  For any given linear function~$w:\R3\to\R{}$,
the subgroup~$G_w\subset\SO(3)$ of rotations that preserve~$w$ is
isomorphic to~$\SO(2)=S^1$.  The induced representation of~$G_w$ 
on~$\cH_3(\R3)\simeq\R7$ is the sum of four 
$\G_w$-irreducible subspaces, $V_0^w$, $V_1^w$, $V_2^w$, and $V_3^w$, 
where $V_0^w$ has dimension~$1$ and is the trivial representation and, 
for~$k>0$, $V_k^w$ has dimension~$2$ and is the representation on
which a rotation by an angle of~$\alpha$ in~$G_w$ acts as a rotation
by an angle of~$k\alpha$.  

For example, the proof of Proposition~\ref{prop: nontrivial stabilizers} 
lists an explicit basis for~$V_k^z$ for~$0\le k\le 3$.  Note that
a cubic~$h\in\cH_3(\R3)$ is linear in~$z$ 
if and only if it lies in~$V_2^z+V_3^z$.  By symmetry, it follows
that a cubic in $\cH_3(\R3)$ is linear in a variable~$w$ 
if and only if it lies in~$V_2^w+V_3^w$.

\begin{proposition}\label{prop: linear cubics}
The set of cubics in~$\cH_3(\R3)$ that are linear in some variable 
is a closed semi-analytic variety of codimension~$1$ in~$\cH_3(\R3)$ 
and consists of the cubics~$h\in\cH_3(\R3)$ 
for which the projective plane curve~$h=0$ has a real singular point.   

Any cubic that is linear in two distinct variables is reducible
and is on the $\SO(3)$-orbit of
$$
h = 3s\,xyz + sp\,(x^3-3xy^2),
$$
which, in addition to being linear in~$z$, is linear in~$w = y+pz$
as well.  When~$s$ and $p$ are nonzero, 
this cubic is not linear in any other variables.

Any cubic that is linear in three distinct variables is 
on the $\SO(3)$-orbit of~$3s\,xyz$ for some~$s\ge0$.
\end{proposition}

\begin{proof}
A cubic~$h$ is linear in a direction~$w$ if and only if the direction
generated by~$w$ is a singular point of the projectivized curve~$h=0$
in~$\bbR\bbP^2$.  Thus, the first statement follows, since the
set of real cubic curves with a real singular point is a semi-analytic
set of codimension~$1$.  If there are
two distinct singular points, then the curve~$h=0$ must be a union
of a line with a conic.  If there are three distinct
singular points, then the curve~$h=0$ must be the union of three
nonconcurrent lines.
Further details are left to the reader.
\end{proof}

\subsection{Continuous symmetry}\label{ssec: so2 symmetry}
I now want to consider those \sLag{} submanifolds ~$L\subset\C3$ 
whose cubic second fundamental form has an~$\SO(2)$ symmetry 
at each point.  

\begin{example}[$\SO(3)$-invariant \sLag{} submanifolds]
\label{ex: HL SO3 slags}
By looking for \sLag{} submanifolds of~$\C3 = \R3 +\iC\,\R3$ 
that are invariant under the `diagonal' action of~$\SO(3)$ on
the two $\R3$-summands, Harvey and Lawson~\cite{MR85i:53058}
found the following examples:
$$
L_c = \bigl\{\,(s+\iC\,t)\ub \mid \ub\in S^2\subset\R3,
            \ t^3{-}3s^2t = c^3\,\bigr\}.
$$
Here,~$c$ is a (real) constant.  Note that~$L_0$ is the union of 
three \sLag{} $3$-planes.  When~$c\not=0$, the submanifold~$L_c$
has three components and each one is smooth and complete.  In fact,
these three components are isometric, as scalar multiplication in~$\C3$
by a nontrivial cube root of unity permutes them cyclically. 
Each of these components is asymptotic to
one pair of $3$-planes drawn from~$L_0$.  The $\SO(3)$-stabilizer 
of a point of~$L_c$ is isomorphic to~$\SO(2)$, so it follows that 
the fundamental cubic at each point has at least an $\SO(2)$-symmetry.
It is not difficult to verify that this cubic is nowhere vanishing
on~$L_c$.  Note also that, for~$\lambda$ real and nonzero,
~$\lambda\cdot L_c = L_{\lambda c}$, so that, up to scaling, 
all of the~$L_c$ with~$c\not=0$ are isometric.
\end{example}

\begin{theorem}\label{thm: so2 symmetry}
If~$L\subset\C3$ is a connected \sLag{} submanifold whose cubic fundamental
form has an $\SO(2)$ symmetry at each point, then either~$L$ is
a $3$-plane or else $L$ is, up to rigid motion, an open subset of
one of the Harvey-Lawson examples.
\end{theorem}

\begin{proof}
Let~$L\subset\C3$ satisfy the hypotheses of the theorem. 
If the fundamental cubic~$C$ vanishes identically, 
then~$L$ is a $3$-plane, so assume that it does not.  
The locus where~$C$ vanishes is a proper real-analytic subset of~$L$,
so its complement~$L^*$ is open and dense in~$L$.  Replace~$L$ by
a component of~$L^*$, so that it can be assumed that~$C$ is nowhere
vanishing on~$L$.%
\footnote{ In principle, this strategy could cause problems, but, as it
is eventually going to be shown that~$L=L^*$ in the general case
anyway, no problems ensue.}

By Proposition~\ref{prop: nontrivial stabilizers}, since the stabilizer
of~$C_x$ is $\SO(2)$ for all~$x\in L$, there is a positive (real-analytic) 
function~$r:L\to \bbR^+$ with the property that the equation
\begin{equation}\label{eq: so2 C normalized}
C = r\,\omega_1\,\bigl(2\,{\omega_1}^2 -3\,{\omega_2}^2-3\,{\omega_3}^2\bigr)
\end{equation}
defines an~$\SO(2)$-subbundle~$F\subset P_L$ 
of the adapted coframe bundle~$P_L\to L$.    
On the subbundle~$F$, the following identities hold:
\begin{equation}\label{eq: beta for so2 C normalized}
\begin{pmatrix}
\beta_{11}&\beta_{12}&\beta_{13}\\
\beta_{21}&\beta_{22}&\beta_{23}\\
\beta_{31}&\beta_{32}&\beta_{33}
\end{pmatrix}
= \begin{pmatrix}
2r\,\omega_1& -r\,\omega_2&-r\,\omega_3\\
-r\,\omega_2&-r\,\omega_1&0\\
-r\,\omega_3&0&-r\,\omega_1
\end{pmatrix}.
\end{equation}

Moreover, because~$F$ is an~$\SO(2)$-bundle, relations of the form
\begin{equation}\label{eq: alpha for so2 C normalized}
\begin{split}
\alpha_{21} &= t_{21}\,\omega_1 +t_{22}\,\omega_2 +t_{23}\,\omega_3 \\
\alpha_{31} &= t_{31}\,\omega_1 +t_{32}\,\omega_2 +t_{33}\,\omega_3
\end{split}
\end{equation}
hold on~$F$ for some functions~$t_{ij}$.  Moreover, 
for~$i=1$, $2$, $3$ there exist functions~$r_i$ on~$F$ so that
\begin{equation}\label{eq: dr for so2 C normalized}
\d r = r_i\,\omega_i\,.
\end{equation}

Substituting the relations~\eqref{eq: beta for so2 C normalized},
~\eqref{eq: alpha for so2 C normalized},
and~\eqref{eq: dr for so2 C normalized} into the identities
\begin{equation}\label{eq: dbeta for so2 C normalized}
\d\beta_{ij} = -\beta_{ik}\w\alpha_{kj}-\alpha_{ik}\w\beta_{kj}\,
\end{equation}
and using the identities~$\d\omega_i = -\alpha_{ij}\w\omega_j$ then
yields polynomial relations among these quantities that can be solved,
leading to relations of the form
\begin{equation} \label{eq: first identities for so2 C normalized}
\begin{split}
\alpha_{21} &= t\,\omega_2\,,\\
\alpha_{31} &= t\,\omega_3\,,\\
\d r &= -4rt\,\omega_1\,,
\end{split}
\end{equation}
where, for brevity, I have written~$t$ for~$t_{22}$.  

Note that \eqref{eq: first identities for so2 C normalized}
implies that~$d\omega_1 = 0$.  Differentiating the
last equation in~\eqref{eq: first identities for so2 C normalized}
implies that there exists a function~$u$ on~$F$ so that
\begin{equation}
\label{eq: second identities for so2 C normalized}
\d t = u\,\omega_1\,.
\end{equation}
Substituting~\eqref{eq: first identities for so2 C normalized}
and~\eqref{eq: second identities for so2 C normalized} into the
identities
\begin{equation}\label{eq: dalpha for so2 C normalized}
\d\alpha_{ij} = -\alpha_{ik}\w\alpha_{kj}+\beta_{ik}\w\beta_{kj}\,
\end{equation}
and expanding, again using the 
identities~$\d\omega_i = -\alpha_{ij}\w\omega_j$, yields the relations
\begin{equation}
\label{eq: third identities for so2 C normalized}
u = (3r^2-t^2),
\qquad\qquad d\alpha_{23} = (t^2+r^2)\,\omega_2\w\omega_3\,.
\end{equation}
Differentiating these last equations yields only identities.

The structure equations found so far can be summarized as follows:
$F\to L$ is an~$\SO(2)$ bundle on which the 1-forms~$\omega_1,\omega_2,
\omega_3, \alpha_{23} (= -\alpha_{32})$ are a basis.  
They satisfy the structure equations
\begin{equation}\label{eq: str eqs for so2 C normalized}
\begin{split}
\d\omega_1 &= 0,\\
\d\omega_2 &= t\,\omega_1\w\omega_2 - \alpha_{23}\w\omega_3\,,\\
\d\omega_3 &= t\,\omega_1\w\omega_3 + \alpha_{23}\w\omega_3\,,\\
\d\alpha_{23}  &= (t^2+r^2)\,\omega_2\w\omega_3\,,\\
\noalign{\vskip3pt}
\d r &= -4rt\,\omega_1\,,\\
\d t &= (3r^2-t^2)\,\omega_1\,,
\end{split}
\end{equation}
and the exterior derivatives of these equations are identities.

These equations imply that 
$\d\bigl(r^{3/2}+r^{-1/2}t^2\bigr) = 0$.  Since~$L$ and~$F$ are 
connected, it follows that there is a constant~$c>0$ so that
$r^{3/2}+r^{-1/2}t^2 = c^{-3/2}$.  Consequently, there is 
a function~$\theta$ that is well-defined on~$L$ that satisfies
$$
r^{3/4} = c^{-3/4}\cos 3\theta,\qquad\qquad
r^{-1/4}t = c^{-3/4}\sin 3\theta.
$$
and the bound~$|\theta|<\pi/6$.  It then follows from the last
two equations of~\eqref{eq: str eqs for so2 C normalized} that
$$
\omega_1 = c\,\frac{\d\theta}{(\cos 3\theta)^{4/3}}\,.
$$
Moreover, setting~$\eta_i = c^{-1}(\cos3\theta)^{1/3}\,\omega_i$ for
$i = 2$ and~$3$ yields
$$
\d\eta_2 = -\alpha_{23}\w\eta_3\,,\qquad
\d\eta_3 =  \alpha_{23}\w\eta_2\,,\qquad
\d\alpha_{23} = \eta_2\w\eta_3\,,
$$
which are the structure equations of the metric of constant curvature~$1$
on~$S^2$.  

Conversely, if~$d\sigma^2$ is the metric of constant curvature~$1$
on~$S^2$, then, on the product~$L = (-\pi/6,\pi/6)\times S^2$, consider
the quadratic and cubic forms defined by
$$
g = c^2\,\frac{{\d\theta}^2 + \cos^23\theta\,d\sigma^2}
                 {(\cos3\theta)^{8/3}}
\qquad\text{and}\qquad
C = c^2\,\frac{{2\,\d\theta}^3 - 3\cos^23\theta\,\d\theta\,d\sigma^2}
                {(\cos3\theta)^{8/3}}.
$$
The metric~$g$ is complete and the pair~$(g,C)$ satisfy
the Gauss and Codazzi equations that ensure that~$(L,g)$ can be 
isometrically embedded as a \sLag{} $3$-fold in~$\C3$ inducing~$C$
as the fundamental cubic.  Thus, for each value of~$c$, there
exists a corresponding \sLag{} $3$-fold that is complete and unique
up to \sLag{} isometries of~$\C3$.  

Since the parameter~$c$ is accounted for by dilation in~$\C3$, 
it now follows that these \sLag{} $3$-folds are the Harvey-Lawson
examples, as desired.  Note that since these are complete and
since~$r$ is nowhere vanishing, it follows that~$L^*=L$ for
the Harvey-Lawson examples, and hence for all examples.
\end{proof}

\subsection{$A_4$ symmetry}\label{ssec: A4 symmetry}

Now consider those \sLag{} submanifolds ~$L\subset\C3$ 
whose fundamental cubic has an~$A_4$-symmetry 
at each point.  Unfortunately, I cannot begin the discussion by
providing a nontrivial example.

\begin{theorem}\label{thm: A4 symmetry}
The only \sLag{} submanifold of~$\C3$ whose fundamental
cubic has an $A_4$-symmetry at each point is a \sLag{} $3$-plane.
\end{theorem}

\begin{remark}
It is interesting to compare the results of Theorems~\ref{thm: so2 symmetry}
and~\ref{thm: A4 symmetry}.  The $\SO(3)$-orbits consisting of
the~$h\in \cH_3(\R3)$ that have an~$\SO(2)$ symmetry 
form a cone of dimension~$3$ in~$\cH_3(\R3)$, while the ones with
an~$A_4$-symmetry form a cone of dimension~$4$ in~$\cH_3(\R3)$.  Thus,
one might expect, based on  `equation counting',
that the the condition of having all cubics have a $\SO(2)$-symmetry
would have fewer solutions than the condition of having all
cubics have an~$A_4$-symmetry.  However, just the opposite is true.
\end{remark}

\begin{proof}
Let~$L\subset\C3$ be a connected \sLag{} submanifold with
the property that its fundamental cubic~$C$ has an~$A_4$-symmetry
at each point.  If $C$ vanishes identically, then~$L$ is
an open subset of a \sLag{} $3$-plane, so assume that it does not.
Let~$L^*\subset L$ be the dense open subset where~$C$ is nonzero.

By Proposition~\ref{prop: nontrivial stabilizers}, since the stabilizer
of~$C_x$ is $A_4$ for all~$x\in L^*$, there is a positive (real-analytic) 
function~$r:L\to \bbR^+$ for which the equation
\begin{equation}\label{eq: A4 C normalized}
C = 6r\,\omega_1\,\omega_2\,\omega_3
\end{equation}
defines an~$A_4$-subbundle~$F\subset P_L$ over~$L^*$
of the adapted coframe bundle~$P_L\to L$.    On~$F$, 
the following identities hold:
\begin{equation}\label{eq: beta for A4 C normalized}
\begin{pmatrix}
\beta_{11}&\beta_{12}&\beta_{13}\\
\beta_{21}&\beta_{22}&\beta_{23}\\
\beta_{31}&\beta_{32}&\beta_{33}
\end{pmatrix}
= \begin{pmatrix}
0& r\,\omega_3&r\,\omega_2\\
r\,\omega_3&0&r\,\omega_1\\
r\,\omega_2&r\,\omega_1&0
\end{pmatrix}.
\end{equation}

Since~$F$ is an~$A_4$-bundle over~$L^*$, there
are relations
\begin{equation}\label{eq: alpha for A4 C normalized}
\begin{split}
\alpha_{23} &= t_{11}\,\omega_1 +t_{12}\,\omega_2 +t_{13}\,\omega_3 \\
\alpha_{31} &= t_{21}\,\omega_1 +t_{22}\,\omega_2 +t_{23}\,\omega_3 \\
\alpha_{12} &= t_{31}\,\omega_1 +t_{32}\,\omega_2 +t_{33}\,\omega_3
\end{split}
\end{equation}
holding on~$F$ for some functions~$t_{ij}$.  Moreover, there
exist functions~$r_i$ for~$i=1,2,3$ on~$F$ so that
\begin{equation}\label{eq: dr for A4 C normalized}
\d r = r_i\,\omega_i\,.
\end{equation}

Substituting the relations~\eqref{eq: beta for A4 C normalized},
~\eqref{eq: alpha for A4 C normalized},
and~\eqref{eq: dr for A4 C normalized} into the identities
\begin{equation}\label{eq: dbeta for A4 C normalized}
\d\beta_{ij} = -\beta_{ik}\w\alpha_{kj}-\alpha_{ik}\w\beta_{kj}\,
\end{equation}
and using the identities~$\d\omega_i = -\alpha_{ij}\w\omega_j$ then
yields polynomial relations among these quantities that can be solved,
leading to relations
\begin{equation}\label{eq: first identities for A4 C normalized}
\alpha_{ij} = 0\,,\qquad\qquad
\d r = 0.
\end{equation}

Substituting~\eqref{eq: beta for A4 C normalized}
and~\eqref{eq: first identities for A4 C normalized}
into the identities
\begin{equation}\label{eq: dalpha for A4 C normalized}
\d\alpha_{ij} = -\alpha_{ik}\w\alpha_{kj}+\beta_{ik}\w\beta_{kj}\,
\end{equation}
yields $r=0$, contrary to hypothesis.
\end{proof}

\subsection{$S_3$ symmetry}\label{ssec: s3 symmetry}

Now consider the \sLag{} submanifolds of~$\C3$ whose
fundamental cubic has an~$S_3$-symmetry at every point.  In contrast
to the case of $A_4$-symmetry, there clearly are nontrivial examples
of this type.  

\begin{example}[Products]\label{ex: line-and-complex curve}
This example is fairly trivial:  Write~$\C3 = \C1\times\C2$ and
look for \sLag{} submanifolds of the form~$L = \R{}\times\Sigma$, where
$\Sigma\subset\C2$ is a surface.  It is not difficult to see
that there is a unique complex structure on~$\C2$ (not the given one!)
with the property that~$L$ is \sLag{} if and only if~$\Sigma$ is 
a complex curve with respect to this structure.

Explicitly, writing~$z_k = x_k + \iC\,y_k$, one sees~$L=\R{}\times\Sigma$
is \sLag{} for~$\Sigma\subset\C2$ if and only if the $2$-forms
$\d x_2\w\d y_2 + \d x_3\w\d y_3$ and~$\d x_2\w\d y_3 + \d y_2\w\d x_3$
each vanish when pulled back to~$\Sigma$.  Since
$$
(\d x_2\w\d y_2 + \d x_3\w\d y_3) + \iC\,(\d x_2\w\d y_3 + \d y_2\w\d x_3)
= (\d x_2 -\iC\,\d x_3) \w (\d y_2 + \iC\,\d y_3),
$$
these $2$-forms vanish on~$\Sigma$ if and only if~$\Sigma$ is  a
complex curve in~$\C2$ endowed with the
complex structure for which~$u = x_2 -\iC\, x_3$ 
and~$v = y_2 + \iC\, y_3$ are holomorphic.

Now, each of these \sLag{} $3$-folds is easily seen to have its
fundamental cubic be expressible as a cubic polynomial in a pair
of~$1$-forms, from which it follows from 
Proposition~\ref{prop: nontrivial stabilizers} that the 
$\SO(3)$-symmetry group of the cubic at each point is either 
everything (if the cubic vanishes at the given point) or else
isomorphic to~$S_3$.  
\end{example}

\begin{example}[Special Lagrangian cones]\label{ex: slag cones}
A more interesting example is to consider the \sLag{} cones.
Suppose that~$\Sigma\subset S^5$ is a (possibly immersed) 
surface with the property that the cone~$C(\Sigma)\subset\C3$ 
is \sLag{}.  Then it is not difficult to show that the fundamental
cubic of~$C(\Sigma)$ has an~$S_3$-stabilizer at those points
where it is not zero.  (This is because the cubic form uses
only two of the directions.)

The necessary and sufficient conditions on~$\Sigma$ that~$C(\Sigma)$
be \sLag{} are easily stated:  
Let~$\ub:S^5\to\C3$ be the inclusion mapping. 
Define a~$1$-form~$\theta$ on~$S^5$ by~$\theta = J\ub\cdot\d\ub$
and define a~$2$-form~$\Psi$ on~$S^5$ by~$\Psi = \ub\,\lhk\,\Im(\Upsilon)$.
Then~$\Sigma\subset S^5$ has the property that~$C(\Sigma)$ is
\sLag{} if and only if~$\theta$ and~$\Psi$ vanish when pulled back 
to~$\Sigma$.  An elementary application of the Cartan-K\"ahler
theorem~\cite{MR92h:58007}
shows that any real-analytic curve~$\gamma\subset S^5$ to
which~$\theta$ pulls back to be zero lies in an irreducible 
real-analytic surface~$\Sigma$ that satisfies these conditions.
Thus, there are many such surfaces.  (In the terminology of
exterior differential systems, these surfaces depend on 
two functions of one variable.)

In addition, many explicit examples of such surfaces are now known.
For example, in~\cite{math.DG/0005164}, a thorough study is done
of the \sLag{} cones that are invariant under a circle action. 
In fact, the differential equation for these surfaces admits
a B\"acklund transformation and can be formulated as an integrable
system.  In principle, the compact torus solutions can be described
explicitly in terms of $\vartheta$-functions via loop group constructions. 
\end{example}

\begin{example}[Twisted \sLag{} cones]\label{ex: twisted slag cones}
The \sLag{} cones can be generalized somewhat, using a construction
found in~\cite[\S4]{MR92k:53112}.  

Again, let~$\xb:\Sigma\to S^5$ be an immersion of a
simply connected surface for which the
cone on~$\xb(\Sigma)$ is \sLag.  Endow~$\Sigma$ with the metric and
orientation that it inherits from this immersion and 
let~$\ast:\Omega^p(\Sigma)\to\Omega^{2-p}(\Sigma)$ be the associated 
 Hodge star operator.  Since~$\Sigma$
is minimal, it follows that
\begin{equation}\label{eq: x is Lap eign}
{\ast}\d\bigl(\ast\d\xb) + 2\xb = 0.
\end{equation}
Now, let~$b:\Sigma\to\bbR$ be any solution to the second order,
linear elliptic equation
\begin{equation}\label{eq: b is Lap eign}
{\ast}\d\bigl(\ast\d b) + 2b = 0.
\end{equation}
(For example, $b$ could be one of the components of~$\xb$.)  
Equations~\eqref{eq: x is Lap eign}
and~\eqref{eq: b is Lap eign} imply that the vector-valued $1$-form
\begin{equation}
\beta =   \xb\,{\ast}\d b - b\,{\ast}\d\xb
\end{equation}
is closed.  Thus, there exists a $\C3$-valued 
function~$\bb:\Sigma\to\C3$ so that~$\d\bb = \beta$.  

Now, consider the immersion~$X:\R{}\times\Sigma\to\C3$ defined by
\begin{equation}
X = \bb + t\,\xb\,.
\end{equation}
Since~$dX = \xb\,(dt+{\ast}\d b) + t\,\d\xb -b\,{\ast}\d\xb$,
it follows that~$X$ immerses~$\R{}\times\Sigma$
as a \sLag{} $3$-fold in~$\C3$, at least away from the locus~$t = b = 0$
in $\R{}\times\Sigma$, where $X$ fails to be an immersion.
Moreover, at those places where the fundamental cubic of this
immersed submanifold is nonzero, it has $S_3$-symmetry.  

It turns out~\cite{MR92k:53112} that the image~$X(\R{}\times\Sigma)$
determines the data~$\xb:\Sigma\to S^5$ and~$b:\Sigma\to\R{}$
up to a replacement of the form~$(\xb,b)\mapsto(-\xb,-b)$, except
in the case that~$\xb(\Sigma)$ lies in a \sLag{} $3$-plane, in 
which case, $X(\R{}\times\Sigma)$ lies in a parallel $3$-plane.

Note that when~$b=0$, the function~$\bb$ is constant, so that
$X(\R{}\times\Sigma)$ is just a translation of the cone on~$\Sigma$. 
Thus, these examples properly generalize the \sLag{} cones.
I will refer to these examples as \emph{twisted \sLag{} cones}. 

As explained in~\cite{MR92k:53112}, this example can be generalized
somewhat by allowing~$\xb:\Sigma\to S^5$ to be a branched
immersion that is an integral manifold of~$\theta$ and~$\Psi$, but
then one must allow~$b$ to have `pole-type' singularities at
the branch points of the immersion~$\xb$.
\end{example}

\begin{theorem}\label{thm: s3 symmetry}
Suppose that~$L\subset\C3$ is a connected \sLag{} $3$-fold with 
the property that its fundamental cubic at each point has 
an~$S_3$-symmetry. Then either~$L$ is congruent to a 
product~$\R{}\times\Sigma$ as in Example~$\ref{ex: line-and-complex curve}$,
or else~$L$ contains a dense open set~$L^*\subset L$ such that
every point of~$L^*$ has a neighborhood that lies in a twisted \sLag{}
cone~$X(\R{}\times\Sigma)$, as in Example~$\ref{ex: twisted slag cones}$.
\end{theorem}

\begin{proof}
Suppose that~$L\subset\C3$ satisfies the hypotheses of the
theorem.  If the fundamental cubic~$C$ vanishes identically
on~$L$, then~$L$ is a $3$-plane and there is nothing to show,
so suppose that~$C\not\equiv0$. Let~$L^\circ\subset L$ 
be the open dense subset where~$C\not=0$.  

The hypothesis that~$C_x$ has~$S_3$-symmetry at every~$x\in L^\circ$
implies that there is a positive function~$s:L^\circ\to\R{}$ and
an $S_3$-subbundle~$F\subset P_L$ over~$L^\circ$ with 
projection~$\xb:F\to L^\circ\subset\C3$
on which the identity
\begin{equation}\label{eq: s3 cubic form}
C = s\,\bigl({\omega_2}^3-3\,\omega_2{\omega_3}^2\bigr)
\end{equation}
holds.  In particular, the second fundamental form of~$L^\circ$
has the form
\begin{equation}\label{eq: s3 2nd fund form}
\text{I\!I} = J\eb_2\ot s({\omega_2}^2-{\omega_3}^2)
             + J\eb_3\ot s(-2\,\omega_2\omega_3),
\end{equation}
where~$\eb_1,\eb_2,\eb_3$ are the vector-valued functions defined by
the moving frame relation~$\d\xb = \eb_1\,\omega_1 + \eb_2\,\omega_2
+\eb_3\,\omega_3$.

It follows from~\eqref{eq: s3 2nd fund form} that~$L\subset\C3$
is an austere submanifold of dimension~$3$.  By Theorem~4.1 
of~\cite{MR92k:53112}, it follows that either~$L^\circ$ is locally 
the product of a line in~$\C3$ with a minimal surface~$\Sigma$ 
in the orthogonal $5$-plane, or else there exists a dense open
subset~$L^\ast\subset L^\circ$ so that every point of~$L^\ast$
has an open neighborhood in $L^\ast$ that lies in a twisted
cone constructed as in Example~\ref{ex: twisted slag cones}
from a minimal immersion~$\xb:\Sigma\to S^5$ 
and an auxiliary function~$b:\Sigma\to\bbR$ 
satisfying~\eqref{eq: b is Lap eign}.  

Since the group of translations and $\SU(3)$-rotations in~$\C3$
acts transitively on the space of lines, it follows that if~$L$
is locally an orthogonal product~$\R{}\times\Sigma$ and is \sLag, 
then, up to translation by a constant,~$\Sigma$ must be a 
complex curve in the complex $2$-plane~$P$ orthogonal to the linear 
factor, where the complex structure on~$P$ is taken to
be as defined in~Example~\ref{ex: line-and-complex curve}.

On the other hand, if $L$ is not locally an orthogonal product
and so is a twisted cone as described above, then one sees 
from the formula for~$dX$ derived in 
Example~\ref{ex: twisted slag cones} that the immersion~$\xb:\Sigma\to S^5$
must not only be minimal, but must have the property that~$\xb^*\theta = 
\xb^*\Psi=0$ as well, as desired.
\end{proof}

\begin{remark}[Singular behavior]
The reader may be annoyed by the apparent need to restrict 
to the open dense subset~$L^\ast\subset L$.  However, there
are subtle singularity issues that seem to require this.  For
more discussion, see the final pages of~\cite{MR92k:53112}.
\end{remark}

\begin{remark}[Austerity]
Theorem~\ref{thm: s3 symmetry} implies that the austere
\sLag{} 3-folds in~$\C3$ are completely described 
by Examples~\ref{ex: line-and-complex curve} 
and~\ref{ex: twisted slag cones}.  
\end{remark}

\begin{remark}[Generality]
The reader knowledgeable about exterior differential systems may
wonder about the generality of the austere \sLag{} $3$-folds in
the sense of Cartan-K\"ahler theory.  While I have avoided this approach 
to the analysis of these examples in this treatment, I should confess
that I first understood the local geometry of these examples by 
doing a Cartan-K\"ahler analysis.  The obvious exterior differential
system that one writes down for these examples is involutive, 
with Cartan characters~$s_1 = 4$ and $s_2 = s_3 = 0$.  The
characteristic variety of the involutive prolongation consists
of two complex conjugate points, each of multiplicity~$2$.
\end{remark}

\subsubsection{Structure equations}\label{sssec: str eqs s3 symmetry}
For use in the next section, I will record here the structure equations
that one derives for systems of this kind.  I will maintain the
notation established in the proof of Theorem~\ref{thm: s3 symmetry}
for~$L^\circ\subset L$, the function~$s$,
and the $S_3$-bundle $\pi:F\to L^\circ$. Thus, the fundamental
cubic factors as
\begin{equation}\label{eq: s3 C factored}
C = s\,\bigl({\omega_2}^3-3\,\omega_2{\omega_3}^2\bigr)
  = s\,\omega_2\,(\omega_2+\sqrt3\,\omega_3)\,(\omega_2-\sqrt3\,\omega_3).
\end{equation}

In particular, the equations
\begin{equation}\label{eq: beta for s3 C normalized}
\begin{pmatrix}
\beta_{11}&\beta_{12}&\beta_{13}\\
\beta_{21}&\beta_{22}&\beta_{23}\\
\beta_{31}&\beta_{32}&\beta_{33}
\end{pmatrix}
= \begin{pmatrix}
0& 0&0\\
0&s\,\omega_2&-s\,\omega_3\\
0&-s\,\omega_3&-s\,\omega_2
\end{pmatrix}
\end{equation}
hold on~$F$.  Moreover, because~$F$ is an~$S_3$-bundle, 
relations of the form
\begin{equation}\label{eq: alpha for s3 C normalized}
\begin{split}
\alpha_{23} &= t_{11}\,\omega_1 +t_{12}\,\omega_2 +t_{13}\,\omega_3\\ 
\alpha_{31} &= t_{21}\,\omega_1 +t_{22}\,\omega_2 +t_{23}\,\omega_3\\
\alpha_{12} &= t_{31}\,\omega_1 +t_{32}\,\omega_2 +t_{33}\,\omega_3
\end{split}
\end{equation}
hold on~$F$ for some functions~$t_{ij}$.  Also, for~$i=1,2,3$ 
there exist functions~$s_i$ on~$F$ so that
\begin{equation}\label{eq: ds for s3 C normalized}
\d s = s_i\,\omega_i\,.
\end{equation}

Substituting the relations~\eqref{eq: beta for s3 C normalized},
~\eqref{eq: alpha for s3 C normalized},
and~\eqref{eq: ds for s3 C normalized} into the identities
\begin{equation}\label{eq: dbeta for s3 C normalized}
\d\beta_{ij} = -\beta_{ik}\w\alpha_{kj}-\alpha_{ik}\w\beta_{kj}\,
\end{equation}
and using the identities~$\d\omega_i = -\alpha_{ij}\w\omega_j$ then
yields polynomial relations among these quantities that can be solved,
leading to relations of the form
\begin{equation}\label{eq: first identities for s3 C normalized}
\begin{split}
\alpha_{23} &= r_2\,\omega_1 - t_2\,\omega_2 - t_3\,\omega_3\,,\\
\alpha_{31} &=           -3r_2\,\omega_2 - 3r_3\,\omega_3\,,\\
\alpha_{12} &= \phantom{-}3r_3\,\omega_2 - 3r_2\,\omega_3\,,\\
\noalign{\vskip2pt}
\d s &= 3s\bigl(r_3\,\omega_1 -t_3\,\omega_2 + t_2\,\omega_3\bigr)\,,
\end{split}
\end{equation}
where I have renamed the covariant derivative variables in the 
solution for simplicity and symmetry of notation.  

It is worth mentioning that not all of these functions on~$F$
are invariant under the action of the group~$S_3$ on the fibers.
The functions~$s$ and $r_2$ are invariant, the function~$r_3$
and the $1$-form~$\omega_1$ are invariant under the odd order
elements of~$S_3$ but switch sign under an element of order~$2$,
and the complex function~$t = (t_2,t_3):F\to\R2$ is $S_3$-equivariant
when~$\R2$ is appropriately identified with the nontrivial 
irreducible representation of dimension~$2$ of~$S_3$. 
Thus, $s$, $r_2$, ${r_3}^2$, and ${t_2}^2+{t_3}^2$
are all well-defined on~$L$, but $r_3$, for example, 
is only well-defined up to a sign.

Because the exterior differential system mentioned above is involutive,
it can be shown that one can prescribe the functions~$t_2$, $t_3$,
$r_2$, and $r_3$ essentially arbitrarily along any curve on which
${\omega_2}^2+{\omega_3}^2$ is nonzero (the curves defined by the
differential equations~$\omega_2 = \omega_3 = 0$ are
characteristic) and generate a solution.  

The functions $r_2$ and $r_3$ vanish identically if and only if $L$ 
is an orthogonal product. Otherwise, $L$ is (locally) a twisted cone. 

The structure equations derived so far imply that
$$
\omega_2\w \d\omega_2 
= {\ts\frac14}\,(\omega_2\pm\sqrt3\,\omega_3)
        \w\d(\omega_2\pm\sqrt3\,\omega_3)
= -2r_2\,\omega_1\w\omega_2\w\omega_3\,,
$$
so it follows that the three linear factors of~$C$ define integrable
$2$-plane fields on $L^\circ$ if and only if~$r_2 \equiv 0$ (and that
if any one of the three is integrable, then so are the other two).

If one considers the differential system with the
additional condition~$r_2\equiv0$,
one sees that it implies the structure 
equation~$dr_3 = 3{r_3}^2\,\omega_1$ and that the reduced system, 
with this condition added, is still involutive, but now with
Cartan characters~$s_1 = 2$ and~$s_2 = s_3 = 0$.  In fact, the
condition~$r_2=0$ characterizes the \sLag{} cones and (under
the additional condition~$r_3=0$) the orthogonal products.

Finally, note that, when~$r_2=0$, the $2$-dimensional leaves
of the $2$-plane field defined by~$\omega_2 = 0$ 
will not lie in $3$-planes unless~$t_3\equiv0$.  Since
the condition~$t_3 = 0$ is not $S_3$-invariant unless~$t_2 = 0$
as well, it follows that, except in the very special 
case~$r_2 = t_3 = t_2 = 0$, at most one of the three foliations
has its leaves lying in $3$-planes.

It is not difficult to show that, up to congruence, there is only one
example that satisfies~$r_2 \equiv t_3 \equiv t_2 \equiv 0$, namely,
the Harvey-Lawson example~$L\subset\C3$ defined in the standard coordinates
by the equations~$|z_1|^2 = |z_2|^2 = |z_3|^2$ and~$\Im(z_1z_2z_3)=0$.
This cone is cut into surfaces by three distinct families of 
Lagrangian planes.  For example, each element of the
circle of Lagrangian planes defined by the relations
$$
z_1 - \eul^{\iC\theta}\overline{z_2} 
= z_3 - \eul^{-2\iC\theta}\overline{z_3} = 0
$$
meets~$L$ in a $2$-dimensional cone.  One gets the other two families
by permuting the coordinates~$z_i$. 
 
In the case that only one of the three linear divisors of~$C$
defines a foliation by surfaces that lie in~$3$-planes, one
can reduce to a $\bbZ_2$-subbundle of~$F$ by imposing the 
condition that~$t_3\equiv0$. Then, by pursuing the calculation 
of the integrability  conditions, 
one finds that the remaining quantities~$r_3$ and~$t_2$ 
must satisfy the equations
\begin{equation}\label{eq: 2nd identities for s3 C normalized}
\begin{split}
\d r_3 &= 3{r_3}^2\,\omega_1\\
\d t_2 & = 3t_2r_3\,\omega_1 + ({t_2}^2+9{r_3}^2-2s^2)\,\omega_3
\end{split}
\end{equation}

Note that if~$r_3$ vanishes anywhere, it vanishes identically.  As
already mentioned, this is the case of a product.  It is not
difficult to show that any connected example of this kind is 
congruent to an open subset of the \sLag{} $3$-fold~$L_c$ defined
by the equations
$$
y_1 = (x_2 - \iC\,x_3)^2 - (y_2 + \iC\,y_3)^2 - c^2 = 0,
$$
where~$c>0$ is a real parameter.  This meets the circle of
Lagrangian planes defined by
$$
y_1 = \cos\theta\,x_2-\sin\theta\,y_2 = \cos\theta\,x_3+\sin\theta\,y_3 = 0
$$
in congruent surfaces that are hyperbolic cylinders.

On the other hand, if~$r_3$ is nonzero, one can reduce the structure
bundle to a parallelization of~$L$ by imposing the conditions~$t_3=0$
and~$r_3>0$, so assume this.  By the structure equations,
the expression~$G = (s^2 + {t_2}^2 + 9\,{r_3}^2) s^{-2/3} {r_3}^{-4/3}$
is constant on~$L$.  Moreover, one easily sees from the structure
equations that the vector field~$X$ that 
satisfies~$\omega_1(X) = \omega_3(X)=0$ and~$\omega_2(X) 
= s^{-1/3} {r_3}^{-2/3}$ is a symmetry vector field of the system
and hence must correspond to an ambient symmetry of the corresponding
solution.  Since this symmetry must fix the vertex of the cone, it
follows that it is a rotation.  Pursuing this observation, it is
not difficult to show that all of these solutions can be described
as follows:  Let~$\lambda_1\ge\lambda_2>0>\lambda_3$ be real 
numbers satisfying~$\lambda_1+\lambda_2+\lambda_3 = 0$, and consider
the $3$-fold~$L_\lambda\subset\C3$ consisting of the points of the form
$$
\left(r_1\,\eul^{\iC(\pi/6+\lambda_1\,t)},
      r_2\,\eul^{\iC(\pi/6+\lambda_2\,t)},
      r_3\,\eul^{\iC(\pi/6+\lambda_3\,t)}
\right)
$$
where~$t$,~$r_1$, $r_2$, and $r_3$ 
are real numbers satisfying
$\lambda_1\,{r_1}^2+\lambda_2\,{r_2}^2+\lambda_3\,{r_3}^2 = 0$.
Then~$L_\lambda$ is a \sLag{} cone with a foliation by 
$2$-dimensional~$3$-plane slices given by~$\d t = 0$. 
(These $3$-plane slices are all congruent and are Euclidean cones.)
Moreover, every $L$ of the type under discussion is
congruent to~$L_\lambda$ for some~$\lambda$.  

Note, by the way, 
that~$L_\lambda$ is not closed unless the ratios of the~$\lambda_i$ 
are rational.  Thus, for `generic'~$\lambda$, the cone~$L_\lambda$
is dense in the $4$-dimensional cone in~$\C3$ defined by the equations
$$
\lambda_1\,|z_1|^2+\lambda_2\,|z_2|^2+\lambda_3\,|z_3|^2 = 
\Re(z_1z_2z_3) = 0.
$$

Part of the significance of these examples will be explained in the
next section.

\subsection{$\bbZ_2$ symmetry}\label{ssec: z2 symmetry}
Now consider a \sLag{} submanifold~$L\subset\C3$
whose fundamental cubic~$C$ has a~$\bbZ_2$-symmetry at each point.
Equivalently, by Proposition~\ref{prop: redcible cubics}, this is
the same as assuming that the fundamental cubic~$C$ is reducible
at each point.

Several nontrivial examples have already been seen:  In fact, if~$C_x$
has a continuous stabilizer at each~$x$ or if~$C_x$ has an~$S_3$-stabilizer
at each~$x$, then Proposition~\ref{prop: nontrivial stabilizers} shows
that~$C_x$ must be reducible at each point.  In the first
case, the examples are classified by Theorem~\ref{thm: so2 symmetry}
and in the second case, the examples are classified 
by~Theorem~\ref{thm: s3 symmetry}.  However, these examples have 
stabilizer groups strictly larger than $\bbZ_2$, so the interesting
question is whether there exist any other examples. 
By Proposition~\ref{prop: nontrivial stabilizers} 
and Theorem~\ref{thm: A4 symmetry}, any such example~$L$ will have
to have the property that the $\SO(3)$-stabilizer of~$C_x$ is
exactly~$\bbZ_2$ for generic~$x\in L$.

Before discussing explicit examples, I will describe a geometrically
interesting condition that forces there to be a $\bbZ_2$-symmetry
of~$C_x$ for all~$x\in L$.

\begin{proposition}\label{prop: foliated implies C reducible}
Let~$L\subset\C3$ be a \sLag{} submanifold that supports
a smooth codimension~$1$ foliation~$\cS$ with the property that each 
$\cS$-leaf~$S\subset L$ lies in a $3$-plane.  Then $C_x$ is
reducible for all~$x\in L$.  In particular, the $\SO(3)$-stabilizer
of~$C_x$ contains an element of order~$2$.
\end{proposition}

\begin{proof}
It suffices to assume that~$L$ is connected, so do this.

If any $\cS$-leaf~$S$ is planar, even locally,
then this plane must be $\omega$-isotropic and Harvey and Lawson's
Theorem~5.5 of~\S{III} in~\cite{MR85i:53058} implies that~$L$ itself
must contain an open subset of a \sLag{} $3$-plane.  By real-analyticity,
it follows that~$L$ itself is planar and hence that~$C_x$ vanishes
identically for all~$x\in L$.  Thus, from now on, I can assume that
none of the $\cS$-leaves are planar and that~$L$ itself is nonplanar.

Choose~$x\in L$ and restrict~$L$ to
a neighborhood~$U$ on which the foliation 
can be expressed a product, i.e.,~$U = X\bigl( (0,1)\times D\bigr)$ 
for some open domain~$D\subset\bbR^2$, and the $\cS$-leaves
in~$U$ are of the form~$X(t,D)$ for~$t\in(0,1)$.
Then, by hypothesis, for each $t\in(0,1)$,
there exists a unique real $3$-plane~$P(t)\subset\C3$ so that
$U\cap P(t) = X(t,D)$, and the surface~$U\cap P(t)$ 
is $\omega$-isotropic.  Since the surface~$U\cap P(t)$ is
non-planar, the plane~$P(t)$ itself must Lagrangian, although it
cannot be \sLag{}, since, otherwise, the uniqueness
aspect of Harvey and Lawson's Theorem~$5.5$ would imply 
that~$U\subset P(t)$, contradicting the assumption that~$L$ is
not planar.  It is not difficult to see that the curve~$t\mapsto P(t)$ 
must be smooth, since the foliation~$\cS$ is assumed to be smooth.

Now, consider the $\SO(2)$-subbundle~$F\subset P_L$ over~$U$ 
with the property that the vector-valued functions~$\eb_2$ and~$\eb_3$
are an oriented basis of the tangent space to the $\cS$-leaves.  Then
$\omega_1$ is well-defined on~$U$ and vanishes when pulled back to
any~$\cS$-leaf.  

Now, the set of Lagrangian planes that contain $\eb_2$ and~$\eb_3$
is the circle of $3$-planes that contain $\eb_2$, $\eb_3$ 
and that are contained in the span of~$\eb_1,\,J\eb_1,\,\eb_2,\,\eb_3$.
In particular, $P(t)$ lies in this plane for each leaf~$X(t,D)\subset U$.
Since each leaf~$\omega_1 = 0$ lies in~$P(t)$, it follows that
the second fundamental form
$$
\text{I\!I} = J\eb_1\ot Q_1 + J\eb_2\ot Q_2 + J\eb_3\ot Q_3
$$
has the property that~$Q_2$ and $Q_3$ must vanish when restricted to
the 2-planes defined by~$\omega_1=0$, i.e., it must be true that $Q_2$
and~$Q_3$ are multiples of~$\omega_1$.  However, by Euler's homogeneity 
relation
$$
C = \omega_1\,Q_1 + \omega_2\,Q_2 + \omega_3\,Q_3\,,
$$
it now follows that~$C$ itself must be a multiple of~$\omega_1$, 
i.e., $C$ is reducible at every point of~$U$, as desired.

Finally, by Proposition~\ref{prop: redcible cubics},
the $\SO(3)$-stabilizer of~$C_x$ must contain an element of
order~$2$ for all~$x\in L$.
\end{proof}

\begin{remark}[Non-integrable factors and non-planar foliations]
It is worth pointing out that there are examples of \sLag{} 
$3$-folds~$L\subset\C3$ for which the fundamental cubic~$C$ is reducible, 
but for which the factors of~$C$ do not define codimension~$1$
foliations of~$L$.  In fact, by the discussion 
in~\ref{sssec: str eqs s3 symmetry}, it follows that, for the generic
\sLag{} $3$-fold~$L$ for which the fundamental cubic~$C$ has an
$\SO(3)$-stabilizer isomorphic to~$S_3$, the cubic~$C$ factors
into three linear factors, no one of which defines an integrable 
$2$-plane field.

Moreover, even in the case where~$r_2\equiv0$
(in which case, $L$ is a cone), so that the three factors are
each integrable, the leaves of the three foliations will not lie in
$3$-planes unless~$t_3\equiv0$, which does not hold for the 
general \sLag{} cone.   
\end{remark}

\begin{example}[Lawlor-Harvey]
This example was first found by Lawlor~\cite{MR89m:49077}, and
was subsequently generalized and extended by 
Harvey~\cite[7.78--9]{MR91e:53056}.  While their results are
valid in all dimensions, I will only discuss the dimension~$3$
case.

They show that, for any compact $2$-dimensional ellipsoid~$E\subset P$ 
where~$P\subset\C3$ is a Lagrangian (but not \sLag{}) 
$3$-plane, the \sLag{} extension~$L$ of~$E$ is 
foliated in codimension~$1$ by a $1$-parameter family of 
$2$-dimensional ellipsoids, each of which lies in a $3$-plane.  
By Proposition~\ref{prop: foliated implies C reducible}, it follows
that the fundamental cubic of the Lawlor-Harvey examples must
be reducible at each point, and thus have a symmetry of order~$2$.

It is not difficult to see that, except when the
ellipsoid is a round $2$-sphere, the Lawlor-Harvey examples are
not special cases of either the $\SO(2)$-symmetry examples
or of the $S_3$-symmetry examples.  Thus, it follows that, at
least at a generic point $x\in L$, the $\SO(3)$-stabilizer of~$C_x$
must be isomorphic to~$\bbZ_2$.
\end{example}

\begin{remark}[Joyce's extension]
Dominic Joyce has informed%
\footnote{ private communication, 3 July 2000}
me that, in fact, the Lawlor-Harvey foliation result continues
to hold when~$E$ is any quadric surface in~$P$, 
not necessarily an ellipsoid, or even a non-singular quadric. 
\end{remark}
  
\begin{theorem}\label{thm: Z2 symmetry}
Suppose that~$L\subset\C3$ is a connected \sLag{} $3$-fold
whose fundamental cubic~$C$ is of~$\bbZ_2$-stabilizer type
on an open dense subset~$L^*\subset L$.  Then~$L^*$ has
a codimension~$1$ foliation~$\cS$ such that each~$\cS$-leaf
lies in a $3$-plane and is, moreover, a quadric surface
in that~$3$-plane.  The space of maximally extended \sLag{} 
$3$-folds of this type is finite dimensional 
and, in fact, coincides with the space of Lawlor-Harvey 
examples, as extended by Joyce.
\end{theorem}

\begin{proof}
By assumption, at a generic point~$x\in L$, 
the $\SO(3)$-stabilizer subgroup of~$C_x$ is isomorphic to~$\bbZ_2$.  
Let~$L^\circ\subset L$ be the open, dense subset where this holds.
Then by Proposition~\ref{prop: nontrivial stabilizers}, there 
exist positive functions~$r,s:L^\circ\to\bbR$ with~$r\not=s$
and a $\bbZ_2$-subbundle~$F\subset P_L$ over~$L^\circ$ on which
the following identity holds:
\begin{equation}\label{eq: Z2 C normalized}
C = r\,\omega_1\,(2{\omega_1}^2 - 3{\omega_2}^2 - 3{\omega_3}^2)
    + 6 s\,\omega_1\omega_2\omega_3\,.
\end{equation}
(Of course,~$\pi:F\to L^\circ$ is a double cover and the reader
can just think of the coframing~$\omega$ as being well-defined on
$L^\circ$ up to the ambiguity of replacing~$\omega_2$ and~$\omega_3$
by~$-\omega_2$ and~$-\omega_3$.)

Consequently, on the subbundle~$F$, the following identities hold:
\begin{equation}\label{eq: beta for Z2 C normalized}
\begin{pmatrix}
\beta_{11}&\beta_{12}&\beta_{13}\\
\beta_{21}&\beta_{22}&\beta_{23}\\
\beta_{31}&\beta_{32}&\beta_{33}
\end{pmatrix}
= \begin{pmatrix}
      2r\,\omega_1      & s\,\omega_3-r\,\omega_2 & s\,\omega_2-r\,\omega_3\\
s\,\omega_3-r\,\omega_2 &            -r\,\omega_1 & s\,\omega_1            \\
s\,\omega_2-r\,\omega_3 & s\,\omega_1             &            -r\,\omega_1
\end{pmatrix}.
\end{equation}

Moreover, because~$F$ is a~$\bbZ_2$-bundle, relations of the form
\begin{equation}\label{eq: alpha for Z2 C normalized}
\begin{split}
\alpha_{23} &= t_{11}\,\omega_1 +t_{12}\,\omega_2 +t_{13}\,\omega_3\\
\alpha_{31} &= t_{21}\,\omega_1 +t_{22}\,\omega_2 +t_{23}\,\omega_3\\
\alpha_{12} &= t_{31}\,\omega_1 +t_{32}\,\omega_2 +t_{33}\,\omega_3
\end{split}
\end{equation}
hold on~$F$ for some functions~$t_{ij}$.  Moreover, for~$i=1,2,3$ there
exist functions~$r_i$ and~$s_i$ on~$F$ so that
\begin{equation}\label{eq: dr ds for Z2 C normalized}
\d r = r_i\,\omega_i\,,\qquad\qquad \d s = s_i\,\omega_i\,.
\end{equation}

Substituting the relations~\eqref{eq: beta for Z2 C normalized},
~\eqref{eq: alpha for Z2 C normalized},
and~\eqref{eq: dr ds for Z2 C normalized} into the identities
\begin{equation}\label{eq: dbeta for Z2 C normalized}
\d\beta_{ij} = -\beta_{ik}\w\alpha_{kj}-\alpha_{ik}\w\beta_{kj}\,
\end{equation}
and using the identities~$\d\omega_i = -\alpha_{ij}\w\omega_j$ then
yields polynomial relations among these quantities that can be solved,%
\footnote{ During the derivation 
of~\eqref{eq: first identities for Z2 C normalized}, 
one uses the assumptions that~$r$, $s$ and $r^2-s^2$ are all nonzero.}
leading to relations of the form
\begin{equation}\label{eq: first identities for Z2 C normalized}
\begin{split}
\d r &= 2(s^2+2r^2)t_1\,\omega_1 +(2rst_3+s^2t_2)\,\omega_2
                                 -(2rst_2+s^2t_3)\,\omega_3\,,\\
\d s &= s\bigl(6rt_1\,\omega_1 +(2st_3+rt_2)\,\omega_2
                               -(2st_2+rt_3)\,\omega_3\bigr)\,,\\
\noalign{\vskip3pt}
\alpha_{23} &= {\ts\frac12}(st_2-rt_3)\,\omega_2 
                +{\ts\frac12}(st_3-rt_2)\,\omega_3\\
\alpha_{31} &= -st_2\,\omega_1 + st_1\,\omega_2 - rt_1\,\omega_3\,,\\
\alpha_{12} &= -st_3\,\omega_1 + rt_1\,\omega_2 - st_1\,\omega_3\,,\\
\end{split}
\end{equation}
where, for brevity, I have introduced the
notation
$$
t_1 = -t_{23}/r\,,\qquad t_2 = -t_{21}/s\,,\qquad t_3 = -t_{31}/s\,.
$$  

Using~\eqref{eq: first identities for Z2 C normalized} to expand out
the identities
$$
0 = \d(\d\omega_1) = \d(\d\omega_2) = \d(\d\omega_3) 
   = \d(\d r) = \d(\d s)
$$
yields relations on the exterior derivatives of~$t_1$, $t_2$, and $t_3$.
These can be expressed by the condition that there 
exist functions~$u_1$, $u_2$, and $u_3$ so that the equations
\begin{equation}\label{eq: second identities for Z2 C normalized}
\begin{split}
\d t_1 &= (s\,u_1-3r-3r^2\,{t_1}^2)\,\omega_1 \\
\d t_2 &= -3t_1(r\,t_2-s\,t_3)\,\omega_1 
           +(u_2-{\ts\frac32}r\,{t_2}^2)\,\omega_2
           +(u_3+{\ts\frac32}s\,{t_2}^2)\,\omega_3\,,\\
\d t_3 &= -3t_1(r\,t_3-s\,t_2)\,\omega_1 
           -(u_3+{\ts\frac32}s\,{t_3}^2)\,\omega_2
           -(u_2-{\ts\frac32}r\,{t_3}^2)\,\omega_3
\end{split}
\end{equation}
hold.  Substituting~\eqref{eq: first identities for Z2 C normalized}
and~\eqref{eq: second identities for Z2 C normalized} into the
identities
\begin{equation}\label{eq: dalpha for Z2 C normalized}
\d\alpha_{ij} = -\alpha_{ik}\w\alpha_{kj}+\beta_{ik}\w\beta_{kj}\,
\end{equation}
and expanding, again using the 
identities~$\d\omega_i = -\alpha_{ij}\w\omega_j$, yields
\begin{equation}\label{eq: u2 u3 out Z2 C normalized}
\begin{split}
u_2 &= {\ts\frac12}\bigl(-2r\,{t_1}^2+r\,{t_2}^2-3s\,t_2t_3+r\,{t_3}^2\bigr)
         - r - s\,u_1\,, \\
u_3 &= {\ts\frac12}\bigl(\phantom{-}2s\,{t_1}^2
               -s\,{t_2}^2+3r\,t_2t_3-s\,{t_3}^2\bigr)
         + s + r\,u_1 \,.
\end{split}
\end{equation}

Finally, expanding out the identities~$\d(\d t_1)=\d(\d t_2)=\d(\d t_3)=0$
shows that they are equivalent to the formula
\begin{equation}\label{eq: du1 Z2 C normalized}
\begin{split}
\d u_1 =&\quad 
-2t_1\bigl(3r\,u_1+s\,(-{t_1}^2+2\,{t_2}^2+2\,{t_3}^2)\bigr)\,\omega_1\\
&\quad - \bigl(u_1(rt_2+st_3)+3(rt_3+st_2)(1+{t_1}^2)\bigr)\,\omega_2\\
&\quad + \bigl(u_1(st_2+rt_3)-3(rt_2+st_s)(1+{t_1}^2)\bigr)\,\omega_3\,.
\end{split}
\end{equation}
The exterior derivative of~\eqref{eq: du1 Z2 C normalized} is
an identity.  

For future use, I record the formulae
\begin{equation}\label{eq: dw Z2 C normalized}
\begin{split}
\d \omega_1 &= -s\,(t_3\,\omega_2-t_2\,\omega_3)\w\omega_1\,, \\
\d \omega_2 &=  t_1\,(r\,\omega_2-s\,\omega_3)\w\omega_1
                +{\ts\frac12}(rt_3-st_2)\,\omega_2\w\omega_3\,,\\
\d \omega_3 &=  t_1\,(r\,\omega_3-s\,\omega_2)\w\omega_1
                +{\ts\frac12}(rt_2-st_3)\,\omega_2\w\omega_3\,.
\end{split}
\end{equation}
which follow from the identities~$\d\omega_i = -\alpha_{ij}\w\omega_j$
coupled with~\eqref{eq: first identities for Z2 C normalized}.

At this point, it is worthwhile taking stock of what has been 
accomplished.  Consider the system of quantities
$$
\omega_1,\,\omega_2,\,\omega_3,\,r,\,s,\,t_1,\,t_2,\,t_3,\,u_1\,.
$$
The formulae \eqref{eq: dw Z2 C normalized}, 
\eqref{eq: first identities for Z2 C normalized}, 
\eqref{eq: second identities for Z2 C normalized}, 
and \eqref{eq: du1 Z2 C normalized} express the exterior derivatives
of these quantities as polynomials in these quantities.  Moreover, 
the relation~$\d(\d q)=0$ for~$q$ any one of these quantities follows
by formal expansion and use of the given exterior derivative
formulae.

By a theorem%
\footnote{ This was originally part of Cartan's general theory of 
intransitive pseudo-groups.  In more recent times, this theorem has
been subsumed into the theory of \emph{Lie algebroids}.  For an
introduction, the reader could try 
the Appendix of~\cite{math.DG/0003099}.}
of \'Elie Cartan, for any six
constants~$\bar r,\,\bar s,\,\bar t_1,\,\bar t_2,\,\bar t_3,\,\bar u_1$,
there exists an open neighborhood~$U$ of~$0\in\bbR^3$ 
that is endowed with three linearly independent~$1$-forms~$\omega_i$ 
and six functions~$r,\,s,\,t_1,\,t_2,\,t_3,\,u_1$ 
that satisfy the equations \eqref{eq: dw Z2 C normalized}, 
\eqref{eq: first identities for Z2 C normalized}, 
\eqref{eq: second identities for Z2 C normalized}, 
and \eqref{eq: du1 Z2 C normalized} and also satisfy
$$
r(0)=\bar r,\ \ s(0)=\bar s,
\ \ t_1(0)=\bar t_1,\ \ t_2(0)=\bar t_2,\ \ t_3(0)=\bar t_3,
\ \ u_1(0)=\bar u_1\,.
$$
Moreover these functions and forms are real-analytic and unique 
in a neighborhood of~$0$, up to a real-analytic local diffeomorphism 
fixing~$0$.  

Now, given such a system~$(\omega,r,s,t,u)$ 
on a simply connected~$3$-manifold~$L$,
one can set~$\eta_i =0$, define~$\alpha_{ij}=-\alpha_{ji}$ 
by the last three equations 
of~\eqref{eq: first identities for Z2 C normalized}, 
define~$\beta_{ij}=\beta_{ji}$ 
by the equations~\eqref{eq: beta for Z2 C normalized}, and see
that the affine structure equations
\begin{equation}
\begin{split}
\d\omega_i &= -\alpha_{ij}\w\omega_j +  \beta_{ij}\w\eta_j\,,\\
\d  \eta_i &= - \beta_{ij}\w\omega_j - \alpha_{ij}\w\eta_j\,,\\
\d\beta_{ij} &= -\beta_{ik}\w\alpha_{kj}-\alpha_{ik}\w\beta_{kj}\,,\\
\d\alpha_{ij} &= -\alpha_{ik}\w\alpha_{kj}+\beta_{ik}\w\beta_{kj}\,
\end{split}
\end{equation}
are identities.  Thus, there is an immersion of~$L$,
unique up to translation and $\SU(3)$-rotation, as a \sLag{}
$3$-manifold in~$\C3$ that induces these structure equations.

In particular, it follows that the space of germs of \sLag{} $3$-manifolds
in~$\C3$ whose fundamental cubics are 
of the form~\eqref{eq: Z2 C normalized} is of dimension~$6$.  Moreover,
any two that agree to order~$4$ at a single point must be equal in a
neighborhood.  It is not difficult to argue from this that the
space one gets by reducing modulo the equivalence relation
defined by analytic continuation is a $3$-dimensional singular
space.

Now, the first of the equations~\eqref{eq: dw Z2 C normalized}
shows that the $2$-plane field~$\omega_1=0$ is integrable, moreover, 
the structure equations found so far imply
\begin{equation}
\d\bigl(\eb_2\w\eb_3\w(J\eb_1-t_1\,\eb_1)\bigr)\equiv 0\mod \omega_1\,.
\end{equation}
In particular, the $3$-plane~$\eb_2\w\eb_3\w(J\eb_1-t_1\,\eb_1)$
is constant along each leaf of~$\omega_1$ and, moreover each such 
leaf lies in an affine $3$-plane parallel to this $3$-plane.
Thus, all of these examples are foliated in codimension~$1$ by
$3$-plane sections.  

Moreover, an examination of the structure equations
shows that the space of congruence classes of such $3$-plane sections 
is of dimension~$3$, the same as the dimension of quadric surfaces in
$3$-space. In fact, using the structure equations, it is not difficult 
to show that these $3$-plane sections are, in fact, quadric surfaces.  
For the sake of brevity, I will not include the details of this
routine calculation here.  

It follows that these \sLag{} $3$-folds all belong
to the class of Lawlor-Harvey examples, as extended by Joyce. 
\end{proof}

\begin{corollary}\label{cor: foliated implies Lawlor-Harvey-Joyce}
Any connected \sLag{} $3$-fold~$L\subset\C3$ that is foliated
in codimension~$1$ by $3$-plane sections is an open subset of
a Lawlor-Harvey-Joyce example.
\end{corollary}

\begin{proof}
By Proposition~\ref{prop: foliated implies C reducible}, 
any such~$L$ must have a reducible fundamental cubic~$C$.  
Thus, the $\SO(3)$-stabilizer of~$C_x$ at each point contains 
a~$\bbZ_2$ and so is either isomorphic to~$\SO(2)$, $S_3$,
or~$\bbZ_2$.  

If this stabilizer is isomorphic to~$\SO(2)$
at a generic point, then Theorem~\ref{thm: so2 symmetry} applies, showing
that~$L$ is a Lawlor-Harvey-Joyce example.

If this stabilizer is isomorphic to~$S_3$ at a generic point,
the discussion at the end of \S\ref{sssec: str eqs s3 symmetry} shows 
that the only such examples that are foliated in codimension~$1$
by~$3$-plane sections have the property that these sections
are necessarily (possibly singular) quadric surfaces, so that
such an~$L$ is, again, a Lawlor-Harvey-Joyce example.

Finally, if the stabilizer is isomorphic to~$\bbZ_2$ at
a generic point, then Theorem~\ref{thm: Z2 symmetry} applies.  
\end{proof}

\begin{remark}[Harvey's characterization]
In his proof of Theorem~7.78 in~\cite{MR91e:53056}, Harvey gives
a characterization of the Lawlor-Harvey examples that is closely
related to Proposition~\ref{prop: foliated implies C reducible}.
What he shows is that any \sLag{} $m$-fold~$L\subset\C{m}$ that
meets a certain concurrent family of Lagrangian $m$-planes in a
codimension~$1$ foliation whose leaves are compact must belong
to the family that they construct.  When~$m=3$, 
Corollary~\ref{cor: foliated implies Lawlor-Harvey-Joyce} is more
general than this, since it makes no assumption about the family
of Lagrangian planes that cut~$L$ to produce the foliation
and makes no assumption about compactness 
(or even completeness) of the leaves.

Of course, one expects that the higher dimensional analog
of Corollary~\ref{cor: foliated implies Lawlor-Harvey-Joyce}
holds, i.e., that any connected \sLag{} $m$-fold~$L\subset\C{m}$ 
that is foliated in codimension~$1$ by $m$-plane sections is
necessarily an open subset of a Lawlor-Harvey-Joyce example.
I have not tried to prove this, but it should be straightforward.
\end{remark}

\subsection{$\bbZ_3$ symmetry}\label{ssec: z3 symmetry}

Now consider those \sLag{} submanifolds~$L\subset\C3$ 
whose cubic second fundamental form has an~$\bbZ_3$-symmetry 
at each point.  

\begin{example}\label{ex: slag fan on min sub of S5}
Let~$\Sigma\subset S^5$ be a surface such that the cone on~$\Sigma$
is \sLag, and consider the $3$-fold
$$
L_\Sigma = \bigl\{\,(s+\iC\,t)\ub \mid \ub\in \Sigma,
            \ t^3{-}3s^2t = c\,\bigr\},
$$
where~$c$ is a (real) constant.  This~$L_\Sigma$ is \sLag.
For example, see~\cite{math.DG/0005164}, where a more general
result for \sLag\ cones in~$\C{n}$ is proved.

Note that~$L_\Sigma$ is diffeomorphic to the disjoint union of 
three copies of~$\bbR\times\Sigma$.  In fact, each `end' of
each component of~$L_\Sigma$ is asymptotic to
the cone on~$\lambda\cdot\Sigma\subset S^5$ 
for some~$\lambda$ satisfying~$\lambda^6=1$.

When~$\Sigma$ is not totally geodesic in~$S^5$ the $\SO(3)$-stabilizer 
of the fundamental cubic at a generic point of point of~$L_\Sigma$ 
is isomorphic to~$\bbZ_3$.
\end{example}

\begin{theorem}\label{thm: Z3 symmetry}
If~$L\subset\C3$ is a connected \sLag{} submanifold whose fundamental
cubic has $\bbZ_3$-symmetry at each point of a dense open subset
of~$L$, then~$L$ contains a dense open set~$L^*$
such that every point of~$L^*$ has an open neighborhood in~$L$ 
that is an open subset of one of the \sLag{} $3$-folds
of Example~$\ref{ex: slag fan on min sub of S5}$.
\end{theorem}

\begin{proof}
Let~$L\subset\C3$ satisfy the hypotheses of the theorem.  
The locus of points~$x\in L$ for which the $\SO(3)$-stabilizer
of~$C$ is larger than~$\bbZ_3$ is a proper real-analytic subset of~$L$,
so its complement~$L^*$ is open and dense in~$L$.  
Thus, I can, without loss of generality, replace~$L$ 
by a component of~$L^*$.  In other words, I can assume
that the $\SO(3)$-stabilizer of~$C_x$ is isomorphic to~$\bbZ_3$
for all~$x\in L$.

By Proposition~\ref{prop: nontrivial stabilizers}, since the stabilizer
of~$C_x$ is $\bbZ_3$ for all~$x\in L$, there are positive (real-analytic) 
functions~$r$ and~$s$ on~$L$ with the property that the equation
\begin{equation}\label{eq: Z3 C normalized}
C = r\,\omega_1\,\bigl(2\,{\omega_1}^2 -3\,{\omega_2}^2-3\,{\omega_3}^2\bigr)
    + s\,\bigl({\omega_2}^3-3\,\omega_2{\omega_3}^2\bigr)
\end{equation}
defines a~$\bbZ_3$-subbundle~$F\subset P_L$ 
of the adapted coframe bundle~$P_L\to L$.    Moreover,
the expression~$s-r\,\sqrt2$ is nowhere vanishing on~$L$.

Now, on the 
subbundle~$F$, the following identities hold:
\begin{equation}\label{eq: beta for Z3 C normalized}
\begin{pmatrix}
\beta_{11}&\beta_{12}&\beta_{13}\\
\beta_{21}&\beta_{22}&\beta_{23}\\
\beta_{31}&\beta_{32}&\beta_{33}
\end{pmatrix}
= \begin{pmatrix}
2r\,\omega_1& -r\,\omega_2&-r\,\omega_3\\
-r\,\omega_2&-r\,\omega_1+s\,\omega_2&-s\,\omega_3\\
-r\,\omega_3&-s\,\omega_3&-r\,\omega_1-s\,\omega_2
\end{pmatrix}.
\end{equation}

Moreover, because~$F$ is a~$\bbZ_3$-bundle, relations of the form
\begin{equation}\label{eq: alpha for Z3 C normalized}
\begin{split}
\alpha_{23} &= t_{11}\,\omega_1 +t_{12}\,\omega_2 +t_{13}\,\omega_3\\
\alpha_{31} &= t_{21}\,\omega_1 +t_{22}\,\omega_2 +t_{23}\,\omega_3\\
\alpha_{12} &= t_{31}\,\omega_1 +t_{32}\,\omega_2 +t_{33}\,\omega_3
\end{split}
\end{equation}
hold on~$F$ for some functions~$t_{ij}$.  Moreover, for~$i=1,2,3$ there
exist functions~$r_i$ and~$s_i$ on~$F$ so that
\begin{equation}\label{eq: dr ds for Z3 C normalized}
\d r = r_i\,\omega_i\,,\qquad\qquad \d s = s_i\,\omega_i\,.
\end{equation}

Substituting the relations~\eqref{eq: beta for Z3 C normalized},
~\eqref{eq: alpha for Z3 C normalized},
and~\eqref{eq: dr ds for Z3 C normalized} into the identities
\begin{equation}\label{eq: dbeta for Z3 C normalized}
\d\beta_{ij} = -\beta_{ik}\w\alpha_{kj}-\alpha_{ik}\w\beta_{kj}\,
\end{equation}
and using the identities~$\d\omega_i = -\alpha_{ij}\w\omega_j$ then
yields polynomial relations among these quantities that can be solved,%
\footnote{ During the derivation 
of~\eqref{eq: first identities for Z3 C normalized}, 
one uses the assumptions that~$r$ and $s$ are nonzero
but not the assumption that~$s{-}r\sqrt2$ is nonzero.}
leading to relations of the form
\begin{equation}\label{eq: first identities for Z3 C normalized}
\begin{split}
\d r &= -4rt_1\,\omega_1\,,\\
\d s &= -s\bigl(t_1\,\omega_1 +3t_3\,\omega_2
                               -3t_2\,\omega_3\bigr)\,,\\
\noalign{\vskip3pt}
\alpha_{23} &=            -t_2\,\omega_2  - t_3\,\omega_3\\
\alpha_{31} &=  \phantom{-}t_1\,\omega_3\,,\\
\alpha_{12} &=            -t_1\,\omega_2\,,\\
\end{split}
\end{equation}
where, for brevity, I have introduced the notation
$$
t_1 = t_{23}\,,\qquad t_2 = -t_{12}\,,\qquad t_3 = -t_{13}\,.
$$  

Using~\eqref{eq: first identities for Z3 C normalized} to expand out
the identities
$$
0 = \d(\d\omega_1) = \d(\d\omega_2) = \d(\d\omega_3) 
   = \d(\d r) = \d(\d s)
$$
and also the identities
$$
\d\alpha_{ij} = -\alpha_{ik}\w\alpha_{kj}+\beta_{ik}\w\beta_{kj}\,
$$
yields relations on the exterior derivatives of~$t_1$, $t_2$, and $t_3$.
When these are solved, one finds that there are functions~$u_2$
and~$u_3$ so that the equations
\begin{equation}\label{eq: second identities for Z3 C normalized}
\begin{split}
\d t_1 &= (3r^2-{t_1}^2)\,\omega_1 \\
\d t_2 &= -t_1t_2\,\omega_1  + u_2\,\omega_2 + (u_3+v)\,\omega_3\,,\\
\d t_3 &= -t_1t_3\,\omega_1  - u_2\,\omega_3 + (u_3-v)\,\omega_2
\end{split}
\end{equation}
hold where
$$
v = s^2-{\ts\frac12}(r^2+{t_1}^2+{t_2}^2+{t_3}^2).
$$

Observe that, if one sets $r=0$ in the current structure equations, 
then these become, up to a trivial change of notation, 
the same structure equations as those for the \sLag{} cones 
discussed in~\S\ref{ssec: s3 symmetry}.  This
is a first hint that these examples must be related to the 
\sLag{} cones.%
\footnote{ Also, if one now computes the Cartan characters of 
the na{\"\i}ve exterior differential system that models these 
structure equations, one finds that~$s_1 = 2$ while $s_2 = s_3 = 0$ 
and that this exterior differential system is involutive.
The characteristic variety is a pair of complex conjugate points,
each of multiplicity~$1$.} 

The next observation is that the structure equations
\begin{equation}\label{eq: dr dt for Z3 C normalized}
\d r = -4rt_1\,\omega_1\,,\qquad\text{and}
\qquad \d t_1 = (3r^2-{t_1}^2)\,\omega_1
\end{equation}
are identical (after replacing~$t_1$ by~$t$) to the last two equations 
of~\eqref{eq: str eqs for so2 C normalized}.  In particular, there
must exist a constant~$c>0$ and a function~$\theta$ on~$L$
satisfying the bound~$|\theta|<\pi/6$ so that
\begin{equation}\label{eq: r t for Z3 C normalized}
r^{3/4} = c^{3/4}\cos 3\theta,\qquad\qquad
r^{-1/4}t_1 = c^{3/4}\sin 3\theta.
\end{equation}
It then follows from~\eqref{eq: dr dt for Z3 C normalized} that
\begin{equation}\label{eq: w1 for Z3 C normalized}
\omega_1 = \frac{\d\theta}{c(\cos 3\theta)^{4/3}}\,.
\end{equation}
By dilation in~$\C3$, one can reduce to the case~$c=1$, so
assume this from now on.  

Consider the following expressions:
\begin{equation}\label{eq: eta for Z3 C normalized}
\begin{split}
p &= r^{-1/4}s,\\
q_2 &= r^{-1/4}t_2,\quad q_3 = r^{-1/4}t_3,\\
v_2 &= r^{-1/2}u_2,\quad v_3 = r^{-1/2}u_3,\\
\eta_2 &= r^{1/4}\omega_2,\quad \eta_3 = r^{1/4}\omega_3\,.
\end{split}
\end{equation}
The structure equations derived above show that
\begin{equation}\label{eq: deta for Z3 C normalized}
\begin{split}
\d\eta_2 & = q_2\,\eta_2\w\eta_3\\
\d\eta_3 & = q_3\,\eta_2\w\eta_3\\
\d p & = -3p\,(q_3\,\eta_2 - q_2\,\eta_3)\\
\d q_2 & = \phantom{-}v_2\,\eta_2 + (v_3 + w)\,\eta_3\\
\d q_3 & =           -v_2\,\eta_3 + (v_3 - w)\,\eta_2
\end{split}
\end{equation}
where~$w = \frac12(1+ {q_2}^2 + {q_3}^2) - p^2$. 

In particular,~$d(p^{1/3}\eta_2) = d(p^{1/3}\eta_3) = 0$.
Let~$x\in L$ be fixed and let~$U\subset L$ be an $x$-neighborhood
on which there exist functions~$y_2$ and $y_3$ vanishing at~$x$
that satisfy~$p^{1/3}\eta_2 = \d y_2$ and~$p^{1/3}\eta_3 = \d y_3$.
Then the functions~$(\theta,y_2, y_3)$ are independent on~$U$ and,
by shrinking~$U$ if necessary, I can assume 
that~$(\theta,y_2,y_3)(U)\subset\R3$ is a product open set of the
form~$I\times D$ where~$I\subset(-\frac\pi6,\frac\pi6)$ is a connected
interval and~$D\subset \bbR^2$ is a disc centered on the origin. 
Of course, the functions~$p$, $q_2$, $q_3$, $v_2$ and~$v_3$ can
be regarded as functions on~$D$, since their differentials are linear
combinations of~$\d y_2$ and $\d y_3$.  In fact, these functions
and forms can now be regarded as defined on the open 
set~$(-\frac\pi6,\frac\pi6)\times D$ by simply reading the formulae
above backwards.  Thus, for example
$$
s = r^{1/4} p = (\cos 3\theta)^{1/4} p
$$
and so forth.  This gives quantities~$\omega_i$, $r$, $s$, $t_i$,
and~$u_i$ that are well-defined on all of~$(-\frac\pi6,\frac\pi6)\times D$
and that satisfy the originally derived structure equations.  It 
follows that there is an immersion of~$(-\frac\pi6,\frac\pi6)\times D$
into~$\C3$ as a \sLag{} $3$-fold that extends~$U$ and pulls back
the constructed forms and quantities to agree with the given ones
on~$U$.  The chief difference is that each of the $\theta$-curves
in~$(-\frac\pi6,\frac\pi6)\times D$ is mapped to a complete curve
in~$\C3$.

Next, observe that the equations
\begin{equation}
\left.
\begin{aligned}
  \d\xb &\equiv \eb_1\,\omega_1\\
\d\eb_1 &\equiv J\eb_1\,(2r\omega_1)\\
\d(J\eb_1) &\equiv -\eb_1\,(2r\omega_1)
\end{aligned}
\quad\right\} \mod \omega_2\,,\omega_3\,,
\end{equation}
which are identical to the corresponding equations
in~\S\ref{ssec: so2 symmetry}, then show that the leaves of the 
curve foliation defined by~$\omega_2 = \omega_3 = 0$ are congruent
to the leaves of the corresponding foliation by the $\eb_1$-curves
in~\S\ref{ssec: so2 symmetry}.  

Finally, note that, setting~$\theta = 0$ (i.e., $t_1= 0$ and~$r=1$ 
in the above structure equations on~$(-\frac\pi6,\frac\pi6)\times D$ 
gives an immersion of~$D$ into~$S^5\subset\C3$ with the property
that the cone on the image~$\Sigma$ is a \sLag{} $3$-fold.  Because the
$\theta$-curves meet this surface orthogonally, it follows easily
that the image of $(-\frac\pi6,\frac\pi6)\times D$ is exactly~$L_\Sigma$
as described in Example~\ref{ex: slag fan on min sub of S5}.  Further
details are left to the reader.
\end{proof}

\subsection{The ruled family}\label{ssec: ruled}
In this last subsection, I am going to consider the generality of
the set of ruled \sLag{} $3$-manifolds.  

Examples of ruled \sLag{} $3$-folds in~$\C3$ were constructed
in Harvey and Lawson's original paper~\cite{MR85i:53058}.  These
included products, \sLag{} cones, and conormal bundles of minimal
surfaces in~$\R3$.  All of these families depend on two functions
of one variable in the sense of exterior differential systems.

Harvey and Lawson also showed in~\cite[Theorems~4.9, 4.13]{MR85i:53058} 
how one could deform the conormal bundle of a minimal surface 
in~$\R3$ according to the data of a harmonic function 
on such a surface and obtain more general ruled \sLag{} $3$-folds.
(Borisenko~\cite{MR94f:53099} later gave a somewhat different 
description of the same family.)  These examples
depend on four functions of one variable in the sense of 
exterior differential systems.

On the other hand, the construction in Example~\ref{ex: twisted slag cones} 
of twisted \sLag{} cones provides another family of examples of 
ruled \sLag{} $3$-folds, again depending on four functions of one
variable in the sense of exterior differential systems.  It is easy
to see that this family is distinct from the family described
in~\cite[Theorems~4.9]{MR85i:53058}. 

In this section, I am going to show that the ruled
\sLag{} $3$-folds depend on \emph{six} functions of one variable in the 
sense of exterior differential systems.  Thus, the `explicit' families
that have been constructed so far are only a small part of the 
complete family.  For a different description of ruled \sLag{} $3$-folds,
one should consult Joyce's recent article~\cite{math.DG/0012060}.

\subsubsection{Almost CR-structures and Levi-flatness}
\label{sssec: almost cr-structures}
For the description I plan to give of the ruled \sLag{} submanifold
of~$\C3$, I will need some facts about a generalized notion of
`pseudo-holomorphic curves'.

Recall that an \emph{almost CR-structure} on a manifold~$M$
is a subbundle~$E\subset TM$ of even dimension equipped with 
a complex structure map~$J:E\to E$.  
The \emph{rank} of the CR-structure is the rank of~$E$ as
a complex bundle and the \emph{codimension} of the CR-structure
is the rank of the quotient bundle~$TM/E$.  A (real) curve~$C\subset M$
is said to be an \emph{$E$-curve} if its tangent line at each point
lies in~$E$.  A (real) surface~$S\subset M$ is said to be 
\emph{$E$-holomorphic} if its tangent plane at each point is 
a complex line in~$E$. (In order to avoid confusion, I will not
adopt the standard practice of calling these surfaces `pseudo-holomorphic
curves', or, indeed, curves of any kind.)

An almost CR-structure~$(E,J)$ will be said to be \emph{Levi-flat}
if, for any~$1$-form~$\rho$ on~$M$ that vanishes on~$E$, the $2$-form
$\d\rho$ vanishes on all the $2$-planes that are complex lines in~$E$.  
Note that Levi-flatness is automatic when the codimension of the CR-structure 
is zero and that Levi-flatness generally has no implications about the
`integrability' of the almost CR-structure to a CR-structure, which is 
a different condition altogether.

\begin{proposition}\label{prop: existence of complex curves}
Let~$(E,J)$ be a real-analytic, Levi-flat, almost CR-structure on~$M$
and let~$C\subset M$ be a real-analytic $E$-curve.  Then there is 
an $E$-holomorphic surface~$S\subset M$ that contains~$C$.  This
surface is locally unique in the sense that, for any two such surfaces~$S_1$
and~$S_2$, the intersection~$S_1\cap S_2$ is also an $E$-holomorphic 
surface that contains~$C$.
\end{proposition}

\begin{proof}
This is a straightforward application 
of the Cartan-K\"ahler Theorem~\cite[Chapter~III]{MR92h:58007}
so I will only give the barest details.  This is a local result, so
it suffices to give a local proof.  

Let~$r$ be the rank of~$(E,J)$
and let~$q$ be its codimension.  For any point~$x\in M$, there
is an open $x$-neighborhood~$U\subset M$ on which there exist 
real-analytic $1$-forms~$\theta_1,\dots,\theta_q$ with real values
and~$\omega_1,\dots,\omega_r$ with complex values with the property
that the equations~$\theta_1 = \dots=\theta_q=0$ define the restriction
of~$E$ to~$U$ and with the property that~$\omega_1,\ldots,\omega_r$
are complex linear on~$E$ and are linearly independent over~$\bbC$
at each point of~$U$.  There are identities of the form
$$
\d\theta_\alpha \equiv K_{\alpha ij}\,\omega_i\w\omega_j
+ L_{\alpha ij}\,\omega_i\w\ov{\omega_j} +
\ov{K_{\alpha ij}}\,\ov{\omega_i}\w\ov{\omega_j}
\mod \theta_1,\dots,\theta_q\,.
$$
The hypothesis of Levi-flatness is simply that the
functions~$L_{\alpha ij}$ all vanish identically.  Under this hypothesis,
the real-analytic exterior differential system~$\cI$ generated algebraically 
by the~$\theta_\alpha$ and the real and imaginary parts of the $2$-forms
$\omega_i\w\omega_j$ is involutive and each of the $1$-dimensional
integral elements is regular and lies in a unique $2$-dimensional 
integral element.  Now apply the Cartan-K\"ahler theorem.
\end{proof}

\subsubsection{Oriented lines}\label{sssec: oriented lines}
Since a ruled $3$-manifold in~$\C3$ can be regarded as a surface
in the space of lines in~$\C3$, it is useful to consider the geometry
of this space.  It is slightly more convenient to consider the 
space~$\L$ of oriented lines in~$\C3$, so I will do this.  

The space~$\L$ is naturally diffeomorphic to the tangent
bundle of~$S^5$. Explicitly, the pair~$(\ub,\vb)\in TS^5$ consisting
of a unit vector~$\ub\in S^5$ and a vector~$\vb\in \ub^\perp$ 
corresponds to the oriented line with oriented direction~$\ub$ 
that passes through~$\vb$.  Naturally, I will regard $\ub:\L\to S^5$
and~$\vb:\L\to\C3$ as vector-valued functions on~$\L$.

Thus, a curve~$\gamma:(a,b)\to\L$ can be written 
as~$\gamma(s) = \bigl(\ub(s),\vb(s)\bigr)$ where the curve~$\ub:(a,b)\to S^5$
and the curve~$\vb:(a,b)\to\C3$ satisfy~$\ub(s)\cdot\vb(s) = 0$ 
for all~$s\in(a,b)$.  Such a curve gives rise 
to a mapping~$\Gamma:(a,b)\times\bbR\to\C3$ by the formula
$$
\Gamma(s,t) = \vb(s) + t\,\ub(s).
$$
Assuming that~$\gamma$ is smooth (resp., real-analytic) then
$\Gamma$ is also smooth (resp., real-analytic) and~$\Gamma$ will
be an immersion except on the locus consisting 
of those~$(s,t)\in(a,b)\times\bbR$ 
where~$\bigl(\vb'(s) + t\,\ub'(s)\bigr)\w\ub(s)= 0$.  On the
locus where it is an immersion, the image of~$\Gamma$ is then
a ruled surface in~$\C3$.

More generally, given any smooth (resp., real-analytic) 
map~$\gamma:P\to\L$ where~$P$
is a smooth (resp., real-analytic) manifold
there is an induced smooth (resp., real-analytic) 
map~$\Gamma:P\times\bbR\to \C3$ defined by the same formula as
above.  With the appropriate `generic' assumptions on~$\gamma$, 
the mapping~$\Gamma$ will be an immersion on some open subset
of~$P\times\bbR$ and its image will be a ruled immersion.

There are two natural differential forms on~$\L$ that are
invariant under the complex isometries of~$\C3$.
These are the pair of $1$-forms
$$
\theta = J\ub\cdot \d\ub,\qquad\text{and}\qquad
\tau = J\ub\cdot \d\vb\,.
$$
It is easy to see that $\theta$ and~$\tau$ are linearly independent, 
so their common kernel~$E\subset T\L$ is a bundle of 
rank~$8$.  The significance of these two 1-forms is revealed
in the following result.

\begin{proposition}
A curve~$\gamma:(a,b)\to\L$ is tangent to~$E$ everywhere
if and only if the corresponding ruled `surface'~$\Gamma:(a,b)\times\bbR
\to\C3$ is $\omega$-isotropic.
\end{proposition}

\begin{proof}
This is immediate from the formulae for~$\Gamma$ and~$\omega$.
\end{proof}

\begin{theorem}\label{thm: CR on Lambda}
There is a complex structure~$J$ on~$E\subset T\L$
with the properties 
\begin{enumerate}
\item $(E,J)$ is a real-analytic, Levi-flat
almost CR-structure on~$\L$ that is invariant under the
complex isometries of~$\C3$.
\item Any ruled \sLag{} $3$-fold~$L$ is locally
the image of the~$\Gamma$ associated to an $E$-holomorphic 
surface~$\gamma:S\to\L$.  When~$L$ is not a $3$-plane,
this local representation is either unique or admits at most
one other such representation.
\item For each $E$-holomorphic surface~$\gamma:S\to\L$,
the corresponding map~$\Gamma:S\times\bbR\to\C3$ is ruled and
a \sLag{} immersion on a dense open subset of~$S\times\bbR$.
\item Any non-planar \sLag{} $3$-fold~$L$ that has two distinct 
rulings is a Lawlor-Harvey-Joyce example for which the $2$-dimensional
$3$-plane sections are quadrics that are doubly ruled.
\end{enumerate} 
\end{theorem}

Before going on to the proof of this result, let me state some
immediate corollaries:

\begin{corollary}
A connected \sLag{} $3$-fold~$L\subset\C{3}$ is ruled if and only
if the set~$\L_L$ of lines that intersect~$L$ in nontrivial 
open intervals {\upshape(}which is an analytic subset 
of~$\L${\upshape)} has
dimension at least~$1$.
\end{corollary}

\begin{proof}  I will only sketch the proof, since the details are
straightforward. First, the easy direction: If~$L$ is ruled, 
then the analytic set~$\L_L$ must have dimension~$2$ at least.  

Conversely, if the dimension of~$\L_L$ is at least~$1$,
then it contains an immersed analytic arc~$\gamma:(a,b)\to\L$, 
which generates a ruled surface~$\Gamma(D)\subset L$ for some
appropriate domain~$D\subset(a,b)\times\bbR$.  
The surface~$\Gamma(D)$ must be $\omega$-isotropic
since~$L$ is Lagrangian.  Thus, the arc~$\gamma$ must be an~$E$-curve. 
By Item~$1$ of Theorem~\ref{thm: CR on Lambda} and
Proposition~\ref{prop: existence of complex curves}, this arc lies
in an $E$-holomorphic surface~$\psi:S\subset\L$.  By Item~$3$
of~Theorem~\ref{thm: CR on Lambda}, there is a dense open 
region~$R\subset S\times\bbR$ so that~$\Psi(R)$ is an immersed ruled 
\sLag{} $3$-fold.  It is not hard to see that this $\Psi(R)$ contains at 
least an open subset of~$\Gamma(D)$.  Since by Harvey and Lawson's
Theorem~5.5, the real-analytic $\omega$-isotropic surface~$\Gamma(D)$
lies in a locally unique \sLag{} $3$-fold, it follows that~$\Psi(R)$
and~$L$ must intersect in an open set.  Thus~$L$ is ruled on an
open set.  By real-analyticity and connectedness, it must be ruled
everywhere.
\end{proof}

\begin{corollary}\label{cor: ruled generality}
The ruled \sLag{} $3$-folds in~$\C3$ depend on six functions of one
variable.  
\end{corollary}

\begin{proof}
Combine Theorem~\ref{thm: CR on Lambda} 
and Proposition~\ref{prop: existence of complex curves}.
\end{proof}

\begin{remark}[The characteristic variety]
The characteristic variety of this system turns out to be
a pair of complex conjugate points, each of multiplicity~$3$.
This is particularly interesting for the following reason:
The condition that the fundamental cubic at each point be
singular is a single equation of second order on the \sLag{}
$3$-fold.  Now, as usual, a Lagrangian manifold can be
written as a gradient graph of a potential function, 
in which case, the special Lagrangian condition is a
single second order elliptic equation for the potential. 
Then the condition that the fundamental cubic be singular
is a single \emph{third} order equation for the potential.
By the general theory, the characteristic variety of 
a system consisting of a single elliptic second order
equation and a single third order equation consists of
at most six points by Bezout's Theorem.  Remarkably, the
`singular cubic' system turns out to have such a `maximal'
characteristic variety and to be involutive.  

This must be quite rare.  In fact, so far, I have been unable 
to find another example of a single pointwise equation on the 
second fundamental form that is involutive and has six points 
in its characteristic variety.
\end{remark}

Now for the proof of Theorem~\ref{thm: CR on Lambda}.

\begin{proof}
First, I will define the almost CR-structure on~$\L$ and
show that it is Levi-flat.  Consider the mapping~$\lambda:F\to\L$ that
sends the coframe~$u:T_x\to\C{3}$ to the oriented line spanned
by~$\eb_1(u)$ that passes through~$x$.  Since the structure equations
give
\begin{equation}
\left.
\begin{aligned}
  \d\xb_{\phantom{1}}
   &\equiv \eb_2\,\omega_{2\phantom{1}}  +\eb_3\,\omega_{3\phantom{1}} 
              +J\eb_1\,\eta_{1\phantom{1}} +J\eb_2\,\eta_{2\phantom{1}} 
               +J\eb_3\,\eta_{3\phantom{1}}\\
\d\eb_1 &\equiv \eb_2\,\alpha_{21}  +\eb_3\,\alpha_{31} 
              +J\eb_1\,\beta_{11}+J\eb_2\,\beta_{21}+J\eb_3\,\beta_{31}
\end{aligned}
\right\} \mod \eb_1\,,
\end{equation}
it follows that the ten $1$-forms that appear on the right-hand
side of this equation are $\lambda$-semibasic and it is evident that
$\lambda^*(\theta) = \beta_{11}$ while~$\lambda^*(\tau)=\eta_1$.
The fibers of~$\lambda$ are cosets of the subgroup of the 
motion group that fixes an oriented line in~$\C3$ and hence are 
diffeomorphic to~$\bbR\times\SU(2)$.  In particular, they are connected.

Define complex-valued~$1$-forms on~$F$ by
$$
\zeta_1 = \omega_2+\iC\,\omega_3\,,\quad
\zeta_2 = \eta_2 - \iC\,\eta_3\,,\quad
\zeta_3 = \alpha_{21} + \iC\,\alpha_{31}\,,\quad
\zeta_4 = \beta_{21} - \iC\,\beta_{31}\,.
$$
These forms are $\lambda$-semibasic and satisfy the equations
$$
\d\zeta_1\equiv\cdots\d\zeta_4\equiv 0 
\mod \beta_{11},\eta_1,\zeta_1,\dots,\zeta_4\,,
$$
while
\begin{equation}\label{eq: Levi-flatness of Lambda}
\left.
\begin{aligned}
  \d\beta_{11} &\equiv \zeta_3\w\zeta_4 + \ov{\zeta_3}\w\ov{\zeta_4}\\
 2\,\d\eta_{1} &\equiv \zeta_1\w\zeta_4 - \zeta_2\w\zeta_3 
       + \ov{\zeta_1}\w\ov{\zeta_4}- \ov{\zeta_2}\w\ov{\zeta_3}\ \\
\end{aligned}
\right\} \mod \beta_{11},\eta_1\,.
\end{equation}
Since the fibers of~$\lambda$ are connected, it follows that
there is a (unique) complex structure~$J:E\to E$ so that the 
complex-valued $1$-forms on~$\L$ that are $\bbC$-linear on~$E$
pull back to be linear combinations 
of~$\beta_{11},\eta_1,\zeta_1,\dots,\zeta_4$.   Moreover, the
equations~\eqref{eq: Levi-flatness of Lambda} imply that the
almost CR-structure~$(E,J)$ is Levi-flat, as promised.  This structure
is clearly real-analytic since it is homogeneous under the action of the
complex isometry group on~$\L$.  (Note also, by the way, that 
the equations~\eqref{eq: Levi-flatness of Lambda} also imply that
this almost CR-structure is not integrable.)  This completes the
proof of Item~$1$.

Now suppose that~$L\subset\C3$ is a ruled \sLag{} $3$-fold that
is not a $3$-plane.  Then, on a dense open set, this ruling can be 
chosen to be real-analytic and smooth.  Consider the subbundle~$F_L$ 
of the adapted frame bundle over~$L$ that has~$\eb_1$ tangent to the
ruling direction.  Thus, the curves in~$L$ defined 
by the differential equations~$\omega_2 = \omega_3=0$ are straight lines
and, of course, $\eb_1$ is tangent to these straight lines.  It follows
that~$\d\eb_1\equiv 0\mod\omega_2,\omega_3$.  (In fact, this is necessary
and sufficient that the $\eb_1$-integral curves be straight lines in~$\C3$.)
Since 
$$
\d\eb_1 = \eb_2\,\alpha_{21}  +\eb_3\,\alpha_{31} 
              +J\eb_1\,\beta_{11}+J\eb_2\,\beta_{21}+J\eb_3\,\beta_{31},
$$
it follows, in particular, that~$\beta_{11}\equiv\beta_{21}\equiv\beta_{31}
\equiv 0 \mod \omega_2,\omega_3$.  Since~$\beta_{ij} = h_{ijk}\,\omega_k$,
it follows from this that~$h_{11j}=0$ for~$j = 1$, $2$, and~$3$.  In
particular, the fundamental cubic
$$
C = h_{ijk}\,\omega_i\omega_j\omega_k
$$
is linear in the direction~$\omega_1$.   
Of course, by Proposition~\ref{prop: linear cubics}, it follows that,
at points where~$C$ is non-zero, it is linear in at most three directions.
Moreover, by Proposition~\ref{prop: linear cubics} 
and Theorem~\ref{thm: A4 symmetry}, there is no
non-planar \sLag{} $3$-fold whose cubic is linear in three directions.
Thus, either~$C$ is linear in exactly two directions on a dense open
set, or else it is linear in exactly one direction on a dense open set.

If $C$ is linear in exactly two directions on a dense open set, then,
again by Proposition~\ref{prop: linear cubics}, it follows that~$C$
is reducible at every point and, on a dense open set, cannot have
an $\SO(3)$-stabilizer isomorphic to~$S_3$, since these are not 
linear in two distinct variables.  It follows that the $\SO(3)$-stabilizer
at a generic point is~$\bbZ_2$, so that, by Theorem~\ref{thm: Z2 symmetry},
$L$ must be one of the Lawlor-Harvey-Joyce examples.  Moreover, the
two linearizing directions, since they represent singular points of
the projectivized cubic curve, must lie on the linear factor of~$C$.
Thus, the two possible rulings must lie in the $2$-dimensional 
slices by $3$-planes.  Of course, this can only happen if the 
quadrics that are these slices are doubly ruled.  Conversely, if the
quadrics that are these slices are doubly ruled, then, obviously, $L$
must be doubly ruled as well.  This establishes Item~$4$ (as well as
the fact that a non-planar \sLag{} $3$-fold cannot be triply ruled).

At any rate, note that~$\beta_{11} = h_{11j}\,\omega_j = 0$ and
that
\begin{equation*}
\begin{aligned}
\beta_{21} &= h_{212}\,\omega_2 + h_{213}\,\omega_3\,,\\
\beta_{31} &= h_{312}\,\omega_2 + h_{313}\,\omega_3
            = h_{213}\,\omega_2 - h_{212}\,\omega_3\,,
\end{aligned}
\end{equation*}
where I have used the symmetry and trace conditions on~$h_{ijk}$
together with the condition~$h_{111}=0$.  It follows that
$$
\bigl(\beta_{21}-\iC\,\beta_{31}\bigr)
\w\bigl(\omega_2+\iC\,\omega_3\bigr) = 0.
$$

There are now two cases to deal with.  Either~$\beta_{21}$ 
and $\beta_{31}$ vanish identically or they do not.

Suppose first that~$\beta_{21}\equiv\beta_{31}\equiv 0$. 
In this case, one can, after restricting to a dense
open set, adapt frames so that the fundamental cubic 
has the form
$$
C = h_{222}\,\bigl({\omega_2}^2 - 3\,\omega_2{\omega_3}^2\bigr),
$$
where~$h_{222}>0$.  In particular, the $\SO(3)$-stabilizer
of~$C$ at the generic point is~$S_3$.   
Set~$s = h_{222}$, so that the notation agrees with the notation
established in~\S\ref{ssec: s3 symmetry}.  Looking back at
the structure equations from that section, one sees that
$$
\alpha_{21} + \iC\,\alpha_{31} 
= -3(r_3 + \iC\,r_2)(\omega_2 + \iC\,\omega_3),
$$
In particular, $(\alpha_{21} + \iC\,\alpha_{31})\w(\omega_2 + \iC\,\omega_3)
= 0$.  Since it has already been established that, in this case,
$$
\beta_{11} = \eta_1 = \beta_{21}+\iC\,\beta_{31} = \eta_2 + \iC\,\eta_3 = 0,
$$
it follows immediately that the natural map from the frame bundle
to~$\L$ that sends a coframe~$u\in F_L$ 
to~$\bigl(\xb(u),\eb_1(u)\bigr)$ maps the coframe bundle into an
$E$-holomorphic surface and that this surface is simply the space
of lines of the ruling.  

Now suppose that~$\beta_{21}$ and $\beta_{31}$ 
do not vanish identically.  Then, by restricting to the dense open set
where they are not simultaneously zero, we can reduce frames to arrange
that~$h_{212}=0$, but that~$h_{312}\not=0$. In fact, there will
exist functions~$r\not=0$, $s$, and $t$ so that 
$$
C = 6r\,\omega_1\omega_2\omega_3 
     +s\bigl({\omega_2}^2 - 3\,\omega_2{\omega_3}^2\bigr)
     +t\bigl(3\,{\omega_2}^2\omega_3-{\omega_3}^2\bigr).
$$
This reduces the frames to a finite ambiguity, but I will not worry
about this, since it does not impose any essential difficulty.  Of course,
$s$ and~$t$ cannot vanish identically by Theorem~\ref{thm: A4 symmetry}.
In particular, on this adapted bundle, the following formulae hold:
$$
\begin{pmatrix}
\beta_{11}&\beta_{12}&\beta_{13}\\
\beta_{21}&\beta_{22}&\beta_{23}\\
\beta_{31}&\beta_{32}&\beta_{33}
\end{pmatrix}
= \begin{pmatrix}
0& r\,\omega_3&r\,\omega_2\\
r\,\omega_3&s\,\omega_2 + t\,\omega_3&r\,\omega_1-s\,\omega_3+t\,\omega_2\\
r\,\omega_2&r\,\omega_1-s\,\omega_3+t\,\omega_2&-s\,\omega_2 - t\,\omega_3
\end{pmatrix}.
$$

Now, there are functions~$p_{ij}$ and~$r_i$, $s_i$, and $t_i$ so that
\begin{align*}
dr &= r_1\,\omega_1+ r_2\,\omega_2 + r_3\,\omega_3\,,\\
ds &= s_1\,\omega_1+ s_2\,\omega_2 + s_3\,\omega_3\,,\\
dt &= t_1\,\omega_1+ t_2\,\omega_2 + t_3\,\omega_3\,,\\
\noalign{\vskip 3pt}
\alpha_{32} &= p_{11}\,\omega_1+ p_{12}\,\omega_2 + p_{13}\,\omega_3\,,\\
\alpha_{13} &= p_{21}\,\omega_1+ p_{22}\,\omega_2 + p_{23}\,\omega_3\,,\\
\alpha_{21} &= p_{31}\,\omega_1+ p_{32}\,\omega_2 + p_{33}\,\omega_3\,.
\end{align*}
Just as in previous cases of moving frame analyses, substituting
these equations into the structure equations for~$\d\beta_{ij}$
yields 15 equations on these 18 quantities.  I will not give the
whole solution, since that is not needed for this argument, 
but will merely note that
these equations imply~$p_{21}=p_{31}=0$ and that $p_{22} = p_{33}$
while $p_{23} + p_{32} = 0$.  In particular, this implies
$$
(\alpha_{21} + \iC\,\alpha_{31})\w(\omega_2 + \iC\,\omega_3) = 0,
$$
just as in the first case.  Moreover, since $\beta_{21}= r\,\omega_3$
and~$\beta_{31} = r\,\omega_2$, it also follows that
$$
(\beta_{21} - \iC\,\beta_{31})\w(\omega_2 + \iC\,\omega_3) = 0.
$$
Since it has already been shown that
$$
\beta_{11} = \eta_1 = \eta_2 + \iC\,\eta_3 = 0,
$$
it follows, once again, that the natural map from the frame bundle
to~$\L$ that sends a coframe~$u\in F_L$ 
to~$\bigl(\xb(u),\eb_1(u)\bigr)$ maps the coframe bundle into an
$E$-holomorphic surface and that this surface is simply the space
of lines of the ruling.  

Thus, it has been shown that any ruled \sLag{} $3$-fold is locally
the $3$-fold generated by an $E$-holomorphic surface in~$\L$.

The only thing left to check is that every $E$-holomorphic surface
in~$\L$ generates a \sLag{} $3$-fold in~$\C3$.  However, given
the analysis already done, this is an elementary exercise in 
the moving frame and can be safely left to the reader.
\end{proof}

\begin{remark}[The relation with \sLag{} cones]
It is not difficult to see that there is a Levi-flat almost CR-structure
of codimension~$1$ on~$S^5$ with the property that its holomorphic
surfaces are exactly the links of \sLag{} cones.

In fact, the mapping~$\eb_1:\L\to S^5$ is an almost CR-mapping
in the obvious sense when~$S^5$ is given this almost CR-structure.
In particular, it follows that any ruled \sLag{} $3$-fold is associated
to a \sLag{} cone that one gets by simply translating all of the
ruling lines so that they pass through one fixed point.  It is in
this sense that all of the ruled \sLag{} $3$-folds in~$\C3$ are
`twisted cones' in some sense.  

In light of this fact, it may be that there is a formula for
ruled \sLag{} $3$-folds that is analogous
to the formula for austere $3$-folds given in~\cite{MR92k:53112}.
I have not yet tried to find this.

On the other hand, this relationship shows that there cannot be
a `Weierstrass formula' for the general ruled \sLag{} like the
formula given by Borisenko for the family that he discovered.
The reason is that such a formula would, at the very least, imply
a Weierstrass formula for the links of \sLag{} cones.  However,
it is easy to show that this exterior differential system is equivalent
to a Monge-Ampere system in $5$-dimensions that, by a theorem of Lie,
does not admit a Weierstrass formula.  Thus, the best that one can
hope for is Weierstrass formulae for special cases.
\end{remark}

\begin{remark}[The generalization to the associative case]
As the reader may know, \sLag{} $3$-folds in~$\C3$ are special 
cases of a more general family of calibrated $3$-folds in~$\R7$,
namely, the \emph{associative} $3$-folds as described \S IV 
of~\cite{MR85i:53058}. 

Regarding~$\R7$ as~$\bbR\times\C3$
and using~$x_0$ as the standard linear coordinate on the~$\bbR$-factor,
the $3$-form
\begin{equation}\label{eq:ass form}
\phi = dx_0\w\bigl({\textstyle\frac{\iC}{2}
(\,\d z_1\w\d\overline{z_1}
   +\d z_2\w\d\overline{z_2}+\d z_3\w\d\overline{z_3}\,)}\bigr)
+\Re(\d z_1 \w \d z_2 \w \d z_3)
\end{equation}
is a calibration on~$\R7$, called the \emph{associative calibration}.
The $3$-folds that it calibrates are said to be \emph{associative}.
The associative $3$-folds that lie in the hyperplane~$\{0\}\times\C3$
are exactly the \sLag{} $3$-folds.  However, there are many more 
associative $3$-folds than \sLag{} $3$-folds, 
since, in particular, Harvey and Lawson prove that every connected 
real-analytic surface~$S\subset\R7$ lies in an (essentially) unique
associative $3$-fold~\cite[\S IV.4, Theorem~4.1]{MR85i:53058}

The subgroup of~$\GL(7,\bbR)$ that
stabilizes~$\phi$ is the compact exceptional group~$\G_2$.   It acts
transitively on the oriented lines in~$\R7$ through the origin, and
the $\G_2$-stabilizer of an oriented line is~$\SU(3)$.  In particular,
the group~$\Gamma$ generated by the translations in~$\R7$ 
and the rotations in~$\G_2$ acts transitively on the space~$\L$ 
of oriented lines in~$\R7$.

It can be shown that there is a unique $\Gamma$-invariant almost complex
structure on~$\L$ with the property that any pseudo-holomorphic 
surface~$S\subset L$ defines a ruled associative
$3$-fold~$\Sigma\subset\R7$ (i.e., the oriented union of the 
oriented lines in~$\R7$ that the points of~$S$ represent) and, 
conversely, that if~$\Sigma\subset\R7$ is a ruled associative $3$-fold
that is not a $3$-plane, then the set of oriented lines~$S\subset\L$
that meet~$\Sigma$ in at least an interval is a pseudo-holomorphic
curve in~$\L$.

Further development of this description allows one to give a 
description of the ruled associative $3$-folds of~$\R7$ that
directly generalizes Joyce's description in~\cite{math.DG/0012060}
of the ruled \sLag{} $3$-folds in~$\C3$.

The details of these results will be reported on elsewhere.
\end{remark}

\bibliographystyle{hamsplain}
\bibliography{SLag,Calibrations,rbryant}
\end{document}